\documentclass[12pt, leqno]{amsart}
\usepackage{enumerate}
\usepackage[dvips]{graphics}
\usepackage{amsmath,amssymb,mathrsfs,geometry}
\usepackage[dvipdfmx]{hyperref}
\usepackage{todonotes}
\usepackage{lipsum}
\newcommand{\doublespace}
   {\addtolength{\baselineskip}{0.15\baselineskip}}

\newtheorem{pdef}{Definition}[section]
\newtheorem{thm}[pdef]{Theorem}
\newtheorem{cor}[pdef]{Corollary}
\newtheorem{lem}[pdef]{Lemma}
\newtheorem{exam}[pdef]{Example}

\newtheorem{remark}[pdef]{Remark}
\theoremstyle{definition}
\newtheorem{cond}[pdef]{Condition}
\newtheorem{prop}[pdef]{Proposition}

\newcounter{equationnumber}

\def\bfconv{\text{\,$\boxplus\boxplus$\,}}
\def\bfconvv{\boxplus\boxplus}
\def\V{V}

\renewcommand{\theequation}{\thesection.\arabic{equation}}
\def\mathletters{
    \addtocounter{equation}{1}
    \edef\@currentlabel{\theequation}
    \setcounter{equationnumber}{\value{equation}}
    \setcounter{equation}{0}
    \edef\theequation{\@currentlabel\noexpand\alph{equation}}
    }

\title{Limit theorems in bi-free probability theory}
\author{Takahiro Hasebe, Hao-Wei Huang and Jiun-Chau Wang}
\date {\today}

\address{Department of Mathematics, Hokkaido University, North 10 West 8, Kita-ku, Sapporo 060-0810, Japan}
\email{thasebe@math.sci.hokudai.ac.jp}

\address{Department of Applied Mathematics, National Sun Yat-sen University,
No. 70, Lienhai Road, Kaohsiung 80424, Taiwan, R.O.C.}
\email{hwhuang@math.nsysu.edu.tw}

\address{Department of Mathematics and Statistics, University of Saskatchewan,
106 Wiggins Road, Saskatoon, Saskatchewan S7N 5E6, Canada}
\email{jcwang@math.usask.ca}

\begin{document}
\maketitle \doublespace \pagestyle{myheadings} \thispagestyle{plain}
\markboth{   }{ }

\begin{abstract}
In this paper additive bi-free convolution is defined for general Borel probability measures, and the limiting distributions for sums of bi-free pairs of selfadjoint commuting random variables in an infinitesimal triangular array are determined. These distributions are characterized by their bi-freely infinite divisibility, and moreover, a transfer principle is established for limit theorems in classical probability theory and Voiculescu's bi-free probability theory. Complete descriptions of bi-free stability and fullness of planar probability distributions are also set down. All these results reveal one important feature about the theory of bi-free probability that it parallels the classical theory perfectly well. The emphasis in the whole work is not on the tool of bi-free combinatorics but only on the analytic machinery.
\end{abstract}


\section{Introduction}  The purpose of this paper is to establish an explicit connection between the families of infinitely divisible laws in classical probability theory and bi-free probability theory. As in free probability theory, it is shown that some classical limit theorem has a nice analogue in the bi-free framework.

Denote by $\mathscr{P}_{\mathbb{R}^d}$ the family of Borel probability measures on $\mathbb{R}^d$. The classical convolution $\mu*\nu$ of $\mu$ and $\nu$ in $\mathscr{P}_{\mathbb{R}^d}$ is the probability distribution of the sum of two independent random vectors whose respective distributions are $\mu$ and $\nu$. In the theory of free probability, freeness and free convolution $\boxplus$
are treated as analogues of classical notion of independence and classical convolution for non-commutative random variables, respectively \cite{VDN92}. The latter theory has a \emph{two-faced} extension, which is invented to study pairs of \emph{left} and \emph{right} random variables (also called left and right \emph{faces}) on a free product of complex Hilbert spaces simultaneously \cite{VoicuBF1}. An independent relation put among these pairs is called \emph{bi-freeness}, which gives rise to \emph{bi-free probability theory}. This fascinating theory has grown rapidly and included several interesting findings based on the foundation of free probability theory.

Since the introduction of bi-free probability by Voiculescu in 2013, combinatorial and analytical approaches have been the main research focuses so far \cite{CNSkoufrains1}\cite{CNSkoufrains2}\cite{GHM}\cite{HW}\cite{Skoufranis1}\cite{Skoufranis2}. Given a two-faced pair $(a,b)$ in a $C^*$-probability space $(\mathcal{A},\varphi)$, its \emph{bi-free partial $R$-transform} is defined through its \emph{Cauchy transform} $G_{(a,b)}(z,w)=\varphi((z-a)^{-1}(w-b)^{-1})$ as
\[R_{(a,b)}(z,w)=zR_a(z)+wR_b(w)+1-\frac{zw}{G_{(a,b)}\big(R_a(z)+1/z,R_b(w)+1/w\big)}\] for complex values $z$ and $w$ in a neighborhood of zero, where $R_a$ and $R_b$ are the usual $R$-transforms of $a$ and $b$, respectively \cite{VoicuBF2}. This transform plays the same role as the usual $R$-transform does in the free case. As in the classical situation, the bi-freeness of $(a,b)$ and $(c,d)$ yields the freeness of the left faces $a,c$ and the freeness of the right ones $b,d$. The reader is referred to \cite{subordination}\cite{GS} for some recent developments.

In the present paper, we continue the previous work \cite{HW} and contribute to the research of bi-free harmonic analysis without any emphasis on Voiculescu's original motivation. To accommodate objects like planar probability distributions or integral representations, it is natural to constraint ourselves to commuting and self-adjoint pairs $(a,b)$ in a certain $C^*$-probability space, i.e. $a,b$ are self-adjoint elements whose commutator $ab-ba=0$. This yields an analytic perspective to bi-free probability theory. More precisely, in such a circumstance the Cauchy transform of $(a,b)$ admits an integral form with the joint distribution of $(a,b)$ as its underlying measure. Any measure in $\mathscr{P}_{\mathbb{R}^2}$ with compact support arises in this manner, i.e. it serves as the joint distribution of such a face. The starting point of relating bi-free probability to the classical situation is that by realizing two given compactly supported $\mu,\nu\in\mathscr{P}_{\mathbb{R}^2}$ as joint distributions of two bi-free and commuting pairs, the bi-free partial $R$-transform linearizes the bi-free convolution $\bfconvv$: for $(z,w)$ near $(0,0)$,
\[R_{\mu\bfconvv\nu}(z,w)=R_\mu(z,w)+R_\nu(z,w).\]

An important concept in the study of limit theorems in probability theory is the \emph{infinite divisibility}. A probability distribution on $\mathbb{R}^d$ is infinitely divisible with respect to a binary operation $\star$ on $\mathscr{P}_{\mathbb{R}^d}$ if it can be expressed as the $\star$-convolution of an arbitrary number of copies of identical measures from $\mathscr{P}_{\mathbb{R}^d}$. When $d=1$, this subject has been thoroughly studied by de Finetti, Kolmogorov, L\'{e}vy and Hin\v{c}in in classical probability \cite{Kolmogorov}\cite{Sato}. The free counterpart is also well studied \cite{BerVoicu93}. The theory of infinitely divisible distributions generalizes (free) central limit theorem as they serve as the limit laws for sums of (freely) independent and identically distributed random variables. (Free) Gaussian and (free) Poisson distributions are typical examples of ($\boxplus$-) $*$-infinitely divisible distributions. Distributions of this kind are determined by their characteristic functions or free $R$-transforms, the so-called L\'{e}vy-Hin\v{c}in type representations. Measures in $\mathscr{P}_{\mathbb{R}^2}$ which are $\bfconvv$-infinitely divisible were first studied in \cite{GHM} in the case when they are compactly supported and considered in the general case in \cite{HW}.

The question under investigation in this paper is to provide the criteria for the weak convergence of the sequence
\begin{equation} \label{BFquestion}
\mu_{n1}\bfconv\mu_{n2}\bfconv\cdots\bfconv\mu_{nk_n}\bfconv\delta_{\mathbf{v}_n},
\end{equation}
where $\{k_n\}_{n=1}^\infty$ is a sequence of strictly increasing positive integers, $\{\mu_{nk}\}_{n\geq1,1\leq k\leq k_n}$ is an \emph{infinitesimal triangular array} of measures in $\mathscr{P}_{\mathbb{R}^2}$ and $\delta_{\mathbf{v}_n}$ denotes the dirac measure at the vector $\mathbf{v}_n$ in $\mathbb{R}^2$. Here the infinitesimal condition of a triangular array $\{\mu_{nk}\}$ in $\mathscr{P}_{\mathbb{R}^d}$ means that for any $\epsilon>0$,
\[\lim_{n\to\infty}\max_{1\leq k\leq k_n}\mu_{nk}\big(\{\mathbf{x}\in\mathbb{R}^d:\|\mathbf{x}\|\geq\epsilon\}\big)=0,\]
where $\|\cdot\|$ is the Euclidean norm on $\mathbb{R}^d$. When $d=1$, there is a one-to-one correspondence between the sets of $*$-infinitely divisible laws and $\boxplus$-infinitely divisible laws. Such a correspondence is characterized by the same parameters, a real number and a positive Borel measure on the real line, in the L\'{e}vy-Hin\v{c}in type formulas \cite{BerPata99}\cite{CG08}. We show that under the hypotheses mentioned above the weak convergence of the sequence in (\ref{BFquestion}) is equivalent to that of the sequence
\[\mu_{n1}*\mu_{n2}*\cdots*\mu_{nk_n}*\delta_{\mathbf{v}_n}.\]
Moreover, the limiting distribution in (\ref{BFquestion}) is shown to be $\bfconvv$-infinitely divisible and its bi-free L\'{e}vy-Hin\v{c}in type representation has the same parameters, a vector, a positive semi-definite matrix and a finite positive Borel measure on the plane, as in the classical result, see Theorem \ref{limitthm1}. The classical limit theory for infinitely divisible laws on the plane has its counterpart in bi-free probability theory.

After setting up some basic tools needed for the investigation and proving the generalization of the bi-free convolution of measures in $\mathscr{P}_{\mathbb{R}^2}$ with compact supports to arbitrary ones in Section 2, another useful bi-free L\'{e}vy-Hin\v{c}in integral representation is provided in Section 3. Section 4 is dedicated to building the parallelism between the families of infinitely divisible laws with respect to classical and bi-free convolutions. The investigation of $\bfconvv$-stable laws on $\mathbb{R}^2$ with their domains of attraction is carried out in Section 5, while the fullness of $\bfconvv$-infinitely divisible laws is set down in Section 6.

\section{Preliminaries} We begin with reviewing some definitions and results in \cite{VoicuBF1}\cite{VoicuBF2}. The \emph{Cauchy transform} of a planar Borel probability measure $\mu$, defined as
\[G_\mu(z,w)=\int_{\mathbb{R}^2}\frac{d\mu(s,t)}{(z-s)(w-t)},\] is an analytic function and satisfies the relation
\[G_\mu(\bar{z},\bar{w})=\overline{G_\mu(z,w)}\] on $(\mathbb{C}\backslash\mathbb{R})^2$.
The underlying measure $\mu$ can be recovered from the transform by the \emph{Stieltjes inversion formula}: the
family $\{\mu_\epsilon\}_{\epsilon>0}$ of probability measures on $\mathbb{R}^2$ defined by
\[d\mu_\epsilon(s,t)=-\frac{1}{2}\Re\big[G_\mu(s+i\epsilon,t+i\epsilon)-G_\mu(s+i\epsilon,t-i\epsilon)\big]\frac{dsdt}{\pi^2}\]
converges to $\mu$ weakly as $\epsilon\to0^+$.

A \emph{truncated cone} in the complex plane is a set of the form
\[\Gamma_{\theta,M}:=\{x+iy\in\mathbb{C}:|x|\leq\theta|y|,\;|y|\geq M\},\;\;\;\;\;\theta,M>0.\] For notational convenience, it will be simply denoted by $\Gamma$ in the sequel if $\theta$ and $M$ are known.
Recall from \cite{BerVoicu93} that free Voiculescu's transform $\phi_\nu$ of a measure $\nu\in\mathscr{P}_\mathbb{R}$ is an analytic function and satisfies the relation $F_\nu(\phi_\nu(z)+z)=z$ on a certain truncated cone $\Gamma$, where $F_\nu$ is the reciprocal of the Cauchy transform of $\nu$,
\[G_\nu(z)=\int_\mathbb{R}\frac{d\nu(s)}{z-s}.\] This transform linearizes the free convolution of probability distributions on $\mathbb{R}$.

For a given $\mu\in\mathscr{P}_{\mathbb{R}^2}$, the probability measures defined as $\mu^{(1)}(B)=\mu(B\times\mathbb{R})$ and $\mu^{(2)}(B)=\mu(\mathbb{R}\times B)$ for Borel sets $B\subset\mathbb{R}$ are called the \emph{marginal laws} of $\mu$. The \emph{bi-free
$\phi$-transform} of $\mu$ defined as
\begin{equation} \label{BFphidef}
\phi_\mu(z,w)=\frac{\phi_{\mu^{(1)}}(z)}{z}+\frac{\phi_{\mu^{(2)}}(w)}{w}
+1-\frac{1}{zwG_\mu\left(F_{\mu^{(1)}}^{-1}(z),F_{\mu^{(2)}}^{-1}(w)\right)}
\end{equation} is an analytic function and satisfies the relation
\begin{equation} \label{extphi}
\phi_\mu(\bar{z},\bar{w})=\overline{\phi_\mu(z,w)}
\end{equation}
on $\Gamma^2$ (notice that the denominator of the last term in (\ref{BFphidef}) never vanishes by shrinking the domain $\Gamma$ if necessary).

It was shown in \cite{VoicuBF2} that the bi-free $\phi$-transform
linearizes the bi-free additive convolution $\bfconvv$ of two planar Borel
probability measures $\mu$ and $\nu$ with compact support:
\begin{equation} \label{BFaddformula}
\phi_{\mu\bfconvv\nu}(z,w)=\phi_\mu(z,w)+\phi_\nu(z,w)
\end{equation} for $(z,w)$ in the common domain of these transforms. The marginal laws of the bi-free convolution of compactly supported
probability measures on $\mathbb{R}^2$ can be expressed in terms of the
free convolution of their marginal laws \cite{GHM}: for $j=1,2$,
\begin{equation} \label{bifreeproj1}
(\mu\bfconv\nu)^{(j)}=\mu^{(j)}\boxplus\nu^{(j)}.
\end{equation} A sequence $\{\nu_n\}_{n=1}^\infty\subset\mathscr{P}_{\mathbb{R}^d}$ is said to \emph{converge weakly} to $\nu\in\mathscr{P}_{\mathbb{R}^d}$, denoted by $\nu_n\Rightarrow\nu$, if
\[\lim_{n\to\infty}\int_{\mathbb{R}^d}f\;d\nu_n=\int_{\mathbb{R}^d}f\;d\nu\]
for any bounded and continuous function $f$ on $\mathbb{R}^d$. For any sequence $\{\mu_n\}_{n=1}^\infty\subset\mathscr{P}_{\mathbb{R}^2}$ converging weakly to $\mu\in\mathscr{P}_{\mathbb{R}^2}$, we have for $j=1,2$,
\begin{equation} \label{marginalWconv}
\mu_n^{(j)}
\Rightarrow\mu^{(j)}
\end{equation}

A family $\mathcal{F}$ of Borel probability measures on $\mathbb{R}^d$ is called
\emph{tight} if
\[\lim_{r\to\infty}\sup_{\mu\in\mathcal{F}}\mu(\{\mathbf{x}\in\mathbb{R}^d:\|\mathbf{x}\|>r\})=0.\]
The correlation of the tightness (or weak
convergence) of Borel probability measures on $\mathbb{R}$ and
$\mathbb{R}^2$ and the convergence properties of their free and bi-free $\phi$-transforms are well known \cite{BerPata99}\cite{HW}. To make the presentation accessible for readers of different backgrounds, some results with their proofs in this direction are provided below. Recall that points $z$ in
$\mathbb{C}\backslash\mathbb{R}$ are said to tend to infinity
\emph{non-tangentially}, which is denoted by $z\to_\sphericalangle\infty$, if $z\to\infty$ with $|\Re z/\Im z|$ uniformly bounded.

\begin{prop} \label{tightG}
A family $\mathcal{F}\subset\mathscr{P}_{\mathbb{R}^2}$ is tight if and only if $zwG_\mu(z,w)-1=o(1)$ uniformly for $\mu\in\mathcal{F}$ as
$z,w\to_\sphericalangle\infty$.
\end{prop}

\begin{pf} First suppose that $\mathcal{F}$ is tight. If $\xi\in\mathbb{C}\backslash\mathbb{R}$ with $|\Re\xi/\Im\xi|$
 bounded by $\theta>0$, then
\[\left|\frac{c}{\xi-c}\right|\leq\sqrt{1+\theta^2},\;\;\;\;\;c\in\mathbb{R},\] which is due to the inequality
$(x-c)^2+x^2/\theta^2\geq c^2/(1+\theta^2)$, $x\in\mathbb{R}$. Applying this
inequality to the decomposition
\[\frac{zw}{(z-s)(w-t)}-1=\frac{s}{z-s}+\frac{t}{w-t}+\frac{st}{(z-s)(w-t)}\] shows the existence of some constant $c_\theta>0$
depending on $\theta$ only so that for $|\Re z/\Im z|,|\Re w/\Im w|\leq\theta$,
\[\left|\frac{zw}{(z-s)(w-t)}-1\right|\leq c_\theta,\;\;\;\;\;(s,t)\in\mathbb{R}^2.\]
Hence for any $\mu\in\mathcal{F}$, we have
\[|zwG_\mu(z,w)-1|\leq\frac{r}{|\Im z|}+\frac{r}{|\Im w|}+\frac{r^2}{|\Im z||\Im w|}+c_\theta\mu(\{\|\mathbf{x}\|>r\}),\] which clearly yields the necessity.

Now we prove the sufficiency. Given any $\epsilon>0$, let $M>1$ be large enough so that
\begin{equation} \label{tightGeq}
\sup_{\mu\in\mathcal{F}}|(iy)(iv)G_\mu(iy,iv)-1|<\epsilon
\end{equation} whenever $|y|,|v|\geq M$. Since the inequality
\[\left|\frac{iv}{(iy-s)(iv-t)}\right|\leq1\] holds for $(s,t)\in\mathbb{R}^2$ and $|y|,|v|\geq M$,
it follows that for all $\mu\in\mathcal{F}$,
\[ivG_\mu(iy,iv)=\int_{\mathbb{R}^2}\frac{iv\;d\mu(s,t)}{(iy-s)(iv-t)}
\to\int_{\mathbb{R}^2}\frac{d\mu(s,t)}{iy-s}=G_{\mu^{(1)}}(iy)\] as
$|v|\to+\infty$ by Dominated Convergence
Theorem. Consequently, $iyG_{\mu^{(1)}}(iy)\to1$ uniformly for $\mu\in\mathcal{F}$ as
$|y|\to+\infty$ by letting $|v|\to+\infty$ in (\ref{tightGeq}). This yields the tightness of the marginal sequence
$\{\mu^{(1)}:\mu\in\mathcal{F}\}$ (cf. \cite[Proposition 5.1]{BerVoicu93}). The tightness of $\{\mu^{(2)}:\mu\in\mathcal{F}\}$ can be obtained in a similar way. Now the sufficiency follows since for any $r>0$,
\begin{equation} \label{eqA}
\mu(\{\|(s,t)\|\geq r\})\leq\mu^{(1)}(\{s:|s|\geq r/\sqrt{2}\})+\mu^{(2)}(\{t:|t|\geq r/\sqrt{2}\}).
\end{equation}
\end{pf} \qed

\begin{prop} \label{nontangential}
Let $\mathcal{F}\subset\mathscr{P}_{\mathbb{R}^2}$ be a tight family. Then $F_{\mu^{(j)}}$ is univalent on some common truncated cone $\Gamma$ in $\mathbb{C}$ with image $F_{\mu^{(j)}}(\Gamma)$ containing some contracted cone $\Gamma_{\alpha,L}$ for every $\mu\in\mathcal{F}$ and $j=1,2$. Moreover, $F_{\mu^{(j)}}^{-1}(\xi)=(1+o(1))\xi$ uniformly for $\mu\in\mathcal{F}$ and $F_{\mu^{(j)}}^{-1}(\xi)\to_\sphericalangle\infty$ as $\xi\to\infty$ with $\xi$ staying in $\Gamma_{\alpha,L}$.
\end{prop}

\begin{pf} The existence of such a truncated cone $\Gamma_{\alpha,L}$ and the asymptotic behaviors can be shown by the statements and techniques of \cite[Proposition 5.4]{BerVoicu93} (see also \cite[Proposition 2.6]{BerPata99}.
\end{pf} \qed

In the following the weak convergence of measures in $\mathscr{P}_{\mathbb{R}^2}$ is translated into the asymptotic properties of their bi-free $\phi$-transforms. Recall that the simplified form $\Gamma^2$ denotes the domain of the bi-free $\phi$-transform of a planar probability measure on which (\ref{extphi}) holds.

\begin{prop} \label{weakconvphi}
Let $\{\mu_n\}_{n=1}^\infty\subset\mathscr{P}_{\mathbb{R}^2}$. Then the following assertions are
equivalent.
\begin{enumerate} [$\qquad(1)$]
\item\label{Wconv2} {The sequence $\{\mu_n\}_{n=1}^\infty$ converges weakly
to a planar Borel probability measure $\mu$.}

\item\label{Gconv2} {Functions in the sequence
$\{\phi_{\mu_n}\}_{n=1}^\infty$ are defined on some fixed domain $\Gamma^2$,
converge uniformly on compact sets of $\Gamma^2$ to a function
$\phi$, and $\phi_{\mu_n}(z,w)=o(1)$ uniformly in $n$ as
$z,w\to\infty$ with $(z,w)\in\Gamma^2$.}

\item\label{Pconv2} {Functions in the sequence
$\{\phi_{\mu_n}\}_{n=1}^\infty$ are defined on some fixed domain $\Gamma^2$, $\lim_{n\to\infty}\phi_{\mu_n}(iy,iv)$ exists for
$(iy,iv)\in\Gamma^2$, and $\phi_{\mu_n}(iy,iv)=o(1)$ uniformly in $n$
as $|y|,|v|\to\infty$.}
\end{enumerate}
Moreover, if \eqref{Wconv2} and \eqref{Gconv2} are satisfied, then
$\phi=\phi_\mu$ in $\Gamma^2$.
\end{prop}

\begin{pf} Throughout the proof, we will use the notations $G_n=G_{\mu_n}$, $G_{jn}=G_{\mu_n^{(j)}}$,
$F_{jn}=F_{\mu_n^{(j)}}$, $F_j=F_{\mu^{(j)}}$,
$\phi_{jn}=\phi_{\mu_n^{(j)}}$, and $\phi_j=\phi_{\mu^{(j)}}$ for all $n$ and for $j=1,2$.

First, suppose $\mu_n\Rightarrow\mu$. According to (\ref{marginalWconv}) and \cite[Proposition 5.7]{BerVoicu93}, there exist $\theta,M>0$ so that every $\phi_{jn}$ is defined on
$\Gamma:=\Gamma_{\theta,M}$, $\phi_{jn}\to\phi_j$ uniformly on compact sets of
$\Gamma$ as $n\to\infty$, and $\phi_{jn}(\xi)=o(|\xi|)$ uniformly in $n$ as $\xi\to\infty$ with $\xi\in\Gamma$ for $j=1,2$. Hence
each $\phi_{\mu_n}$ is defined on $\Gamma^2$. On the other hand, the integrands in $G_n$ are uniformly bounded functions of $(s,t)$ for points $(z,w)$ lying in compact sets of
$(\mathbb{C}\backslash\mathbb{R})^2$. This yields the normality of $\{G_n\}$ by Montel's theorem in complex analysis
of several variables. Hence $G_n(F_{1n}^{-1},F_{2n}^{-1})\to G_\mu(F_1^{-1},F_2^{-1})$
uniformly on compact sets of $\Gamma^2$ as $n\to\infty$ and $\phi=\phi_\mu$. To
finish the proof of $\eqref{Wconv2}\Rightarrow\eqref{Gconv2}$, it remains to show that
\[\frac{zwG_n\left(F_{1n}^{-1}(z),F_{2n}^{-1}(w)\right)-1}{zwG_n\left(F_{1n}^{-1}(z),F_{2n}^{-1}(w)\right)}=o(1)\] uniformly in $n$ as $z,w\to\infty$ with
$(z,w)\in\Gamma^2$, which is equivalent to showing the uniform convergence of
$zwG_n(F_{1n}^{-1}(z),F_{2n}^{-1}(w))-1=o(1)$. Note that this can be obtained by applying Proposition \ref{tightG} and Proposition \ref{nontangential} to the identity
\begin{equation} \label{uniformeq}
zwG_n\left(F_{1n}^{-1}(z),F_{2n}^{-1}(w)\right)=F_{1n}^{-1}(z)F_{2n}^{-1}(w)
G_n\left(F_{1n}^{-1}(z),F_{2n}^{-1}(w)\right)\cdot\frac{z}{F_{1n}^{-1}(z)}\cdot\frac{w}{F_{2n}^{-1}(w)}.
\end{equation}

Clearly, \eqref{Gconv2} implies \eqref{Pconv2}. To show $\eqref{Pconv2}\Rightarrow\eqref{Wconv2}$, it suffices
to verify that the sequence $\{\mu_n\}_{n=1}^\infty$ is tight by the established result $\eqref{Wconv2}\Rightarrow\eqref{Gconv2}$. To conclude the proof,
observe that
\begin{align*}
\phi_{\mu_n}(iy,iv)&=\frac{\phi_{1n}(iy)}{iy}+\frac{\phi_{2n}(iv)}{iv}+1-\frac{1}{(iy)\frac{iv}{F_{2n}^{-1}(iv)}F_{2n}^{-1}(iv)G_n
\left(F_{1n}^{-1}(iy),F_{2n}^{-1}(iv)\right)} \\
&\to\frac{\phi_{1n}(iy)}{iy}+1-\frac{1}{(iy)G_{1n}\left(F_{1n}^{-1}(iy)\right)}\;\;\;\;\;\mathrm{as}
\;\;\;|v|\to\infty \\
&=\frac{\phi_{1n}(iy)}{iy},
\end{align*} where Proposition \ref{nontangential} and the fact that $wG_n(z,w)\to G_{1n}(z)$ for
any $z\in\mathbb{C}\backslash\mathbb{R}$ as $w\to\infty$ non-tangentially are used in the limit. This implies that
$\phi_{1n}(iy)=o(|y|)$ uniformly in $n$ as $|y|\to\infty$. Similarly, $\{\phi_{2n}\}$ has the same asymptotic
property. Hence $(iy)(iv)G_n(F_{1n}^{-1}(iy),F_{2n}^{-1}(iv))-1=o(1)$
uniformly in $n$ as $|y|,|v|\to\infty$ by (\ref{uniformeq}), which yields the tightness of $\{\mu_n\}$ by Proposition
\ref{tightG}. The proof is complete.
\end{pf} \qed

We can now define the bi-free convolution of arbitrary planar Borel probability measures $\mu$ and $\nu$.
Choose two sequences $\{\mu_n\}_{n=1}^\infty$ and
$\{\nu_n\}_{n=1}^\infty$ of planar Borel probability measures with
compact support converging to $\mu$ and $\nu$ weakly, respectively. Then Proposition \ref{weakconvphi} shows that $\mu_n\bfconv\nu_n$ weakly converges to a probability measure $\rho$ on $\mathbb{R}^2$ which satisfies
\begin{equation} \label{BFadditive2}
\phi_\rho(z,w)=\phi_\mu(z,w)+\phi_\nu(z,w)
\end{equation} on some $\Gamma^2$. We further deduce from \cite[Proposition 2.5]{HW} the uniqueness of $\rho$.
These discussions lead into the following definition:

\begin{pdef} {\normalfont For any $\mu,\nu\in\mathscr{P}_{\mathbb{R}^2}$, the unique probability measure $\rho$ satisfying the additive identity in (\ref{BFadditive2}) is called the \emph{bi-free additive convolution} of $\mu$ and $\nu$, and is also
denoted by $\mu\bfconv\nu$.}
\end{pdef}

\begin{remark}
\emph{Our approach to the generalization of bi-free additive convolution is based on analytic tools. An (unbounded) operator model for the bi-free convolution of arbitrary probability measures on $\mathbb{R}^2$ is unknown.}
\end{remark}

We can also show by Proposition \ref{weakconvphi} that the operation of bi-free convolution is weakly continuous, namely, if $\{\mu_n\}_{n=1}^\infty$ and $\{\nu_n\}_{n=1}^\infty$ are in $\mathscr{P}_{\mathbb{R}^2}$ weakly converging to $\mu$ and $\nu$, respectively, then $\mu_n\bfconv\nu_n$ weakly converges to $\mu\bfconv\nu$. Finally, we generalize the result in (\ref{bifreeproj1}) to arbitrary measures in $\mathscr{P}_{\mathbb{R}^2}$.

\begin{prop} \label{bifreeproj2}
Let $\mu,\nu\in\mathscr{P}_{\mathbb{R}^2}$. Then for $j=1,2$,
\[(\mu\bfconv\nu)^{(j)}=\mu^{(j)}\boxplus\nu^{(j)}.\]
\end{prop}

\begin{pf}
Choose two sequences $\{\mu_n\}_{n=1}^\infty$ and $\{\nu_n\}_{n=1}^\infty$ of planar Borel probability measures with
compact support converging to $\mu$ and $\nu$ weakly, respectively. Since the projections onto marginals and bi-free convolution are weakly continuous, the result (\ref{bifreeproj1}) passes to the conclusion.
\end{pf}\qed

\section{Bi-free infinite divisibility and L\'{e}vy-Hin\v{c}in representation}
Throughout the remaining part of the paper, points $(s,t)$ in $\mathbb{R}^2$ will be denoted by the bold letter $\mathbf{x}$ and the origin $(0,0)$ will be written as $\mathbf{0}$. We will also denote by the real numbers $\mathbf{v}^{(1)}$ and $\mathbf{v}^{(2)}$
the $s$- and $t$-coordinate of a given vector $\mathbf{v}\in\mathbb{R}^2$.

In classical probability theory, a Borel probability measure $\mu$ on $\mathbb{R}^2$ is $\ast$-infinitely divisible if and only if its characteristic function is of the form (called the \emph{L\'{e}vy-Hin\v{c}in representation})
\begin{equation}\label{LK1}
\widehat{\mu}(\mathbf{u})=\exp\left[i\langle \mathbf{u},\mathbf{v}\rangle-\frac{1}{2}\langle\mathbf{A}\mathbf{u},\mathbf{u}\rangle+ \int_{\mathbb{R}^2}\left(e^{i\langle\mathbf{u},\mathbf{x}\rangle}-1-\frac{i\langle\mathbf{u},\mathbf{x}\rangle}{1+\|\mathbf{x}\|^2}\right) d\tau(\mathbf{x})\right]
\end{equation}
for some vector $ \mathbf{v} \in \mathbb{R}^2$, real positive semi-definite matrix $\mathbf{A}$ and some positive Borel measure $\tau$ on $\mathbb{R}^2$ with the properties that
$\tau(\{\mathbf{0}\})=0$ and $1\wedge\|\mathbf{x}\|^2 \in L^1(\tau)$, where $1\wedge\|\mathbf{x}\|^2:=\min\{1,\|\mathbf{x}\|^2\}$. Conversely, such a triplet $(\mathbf{v},\mathbf{A},\tau)$ generates a probability measure $\mu$ for which \eqref{LK1} holds. The triplet $(\mathbf{v},\mathbf{A},\tau)$ in the representation is unique and called the (classical) \emph{characteristic triplet} of $\mu$, while the measure $\tau$ is called the (classical) \emph{L\'evy measure} of $\mu$. In this case $\mu$ is denoted by $\mu_\ast^{(\mathbf{v},\mathbf{A},\tau)}$. The reader is referred to \cite{limitthm}\cite{Sato} for more details.

Recall that a measure $\mu\in\mathscr{P}_{\mathbb{R}^2}$ is said to be \emph{bi-freely infinitely divisible} if for any $n\in\mathbb{N}$,
it can be expressed as an $n$-fold bi-free convolution of some $\mu_n\in\mathscr{P}_{\mathbb{R}^2}$:
\[\mu=\underbrace{\mu_n\bfconv\cdots\bfconv\mu_n}_{n\;\;\mathrm{terms}}:=\mu_n^{\bfconvv n}.\] Such a measure is
characterized in terms of the functional properties of its bi-free $\phi$-transform \cite[Theorem 4.3]{HW}:
$\mu$ is $\bfconvv$-infinitely divisible if and only if
$\phi_\mu$ extends analytically to
$(\mathbb{C}\backslash\mathbb{R})^2$ and admits an integral representation of the form
\[\phi_\mu(z,w)=\frac{1}{z}\left(\gamma_1+\int_{\mathbb{R}^2}\frac{1+zs}{z-s}\;d\sigma_1(\mathbf{x})\right)+
\frac{1}{w}\left(\gamma_2+\int_{\mathbb{R}^2}\frac{1+wt}{w-t}\;d\sigma_2(\mathbf{x})\right)+\widetilde{D}(z,w),\] where
$(\gamma_1,\gamma_2)\in\mathbb{R}^2$, $\sigma_j$ is a finite, positive Borel measure on $\mathbb{R}^2$ for $j=1,2$,
and
\[\widetilde{D}(z,w)=\int_{\mathbb{R}^2}\frac{\sqrt{1+s^2}\sqrt{1+t^2}}{(z-s)(w-t)}\;d\widetilde{\sigma}(\mathbf{x})\] for some finite Borel signed measure
$\widetilde{\sigma}$ on $\mathbb{R}^2$ satisfying the relations
\begin{equation} \label{sigmadef}
\left\{
\begin{array}{ll}
\frac{t}{\sqrt{1+t^2}}\;\sigma_1=\frac{s}{\sqrt{1+s^2}}\;\widetilde{\sigma}, \\\\
\frac{s}{\sqrt{1+s^2}}\;\sigma_2=\frac{t}{\sqrt{1+t^2}}\;\widetilde{\sigma}, \\\\
\widetilde{\sigma}(\{\mathbf{0}\})^2\leq\sigma_1(\{\mathbf{0}\})\sigma_2(\{\mathbf{0}\}).
\end{array}\right.
\end{equation}
These parameters $\gamma_1$, $\gamma_2$, $\sigma_1$, $\sigma_2$ and $\widetilde{\sigma}$ appearing in the representation are unique. Notice that the first two relations in (\ref{sigmadef}) indicate that
$a:=\sigma_1(\{0\}\times\mathbb{R})=\sigma_1(\{\mathbf{0}\})$,
$b:=\sigma_2(\mathbb{R}\times\{0\})=\sigma_2(\{\mathbf{0}\})$ and $c:=\widetilde{\sigma}(\{st=0\})=
\widetilde{\sigma}(\{\mathbf{0}\})$. An application of Cauchy-Schwarz inequality gives $c^2\leq ab$, i.e. the positive semi-definiteness of the matrix
\begin{equation} \label{A1}
{\normalfont \bf A}=\left(
\begin{array}{cc}
a & c \\
c & b \\
\end{array}
\right).
\end{equation}

In order to get more insights into $\bfconvv$-infinitely divisible distributions, we will derive another integral representation for them. The  representing measure is no longer required to be finite but only positive. First of all, define the positive measure $\tau$ on $\mathbb{R}^2$ as
\begin{equation} \label{taudef}
\tau=\left\{
\begin{array}{ll}
\frac{1+s^2}{s^2}\;\sigma_1 & \hbox{on\;\;\;$\{(s,t)\in\mathbb{R}^2:s\neq0\}$;} \\\\
\frac{1+t^2}{t^2}\;\sigma_2 & \hbox{on\;\;\;$\{(s,t)\in\mathbb{R}^2:t\neq0\}$,}
\end{array}
\right.
\end{equation}
and $\tau(\{\mathbf{0}\})=0$. The relations among $\sigma_1$, $\sigma_2$ and $\widetilde{\sigma}$ in (\ref{sigmadef}) clearly show that $\tau$ is well-defined. Moreover, the
restriction of $\tau$ to the set $\{(s,t)\in\mathbb{R}^2:st\neq0\}$ is equal to
\[\frac{\sqrt{1+s^2}\sqrt{1+t^2}}{st}\chi_{\{st\neq0\}}(s,t)\;\widetilde{\sigma}.\] It is also easy to verify that the function $1\wedge\|\mathbf{x}\|^2$ belongs to $L^1(\tau)$. After these setups, we can rewrite $\phi_\mu$ as
\begin{equation} \label{BFLK}
\phi_\mu(z,w)=\frac{\mathbf{v}^{(1)}}{z}+\frac{\mathbf{v}^{(2)}}{w}+\left(\frac{a}{z^2}
+\frac{c}{zw}+\frac{b}{w^2}\right)+\mathcal{P}(z,w),
\end{equation}
where $\mathbf{v}\in\mathbb{R}^2$ and
\[\mathcal{P}(z,w)=\int_{\mathbb{R}^2}\left[\frac{zw}{(z-s)(w-t)}
-1-\frac{sz^{-1}+tw^{-1}}{1+s^2+t^2}\right]d\tau(s,t).\]
Indeed, by means of the identities
\begin{equation} \label{decom}
\frac{zs}{z-s}=\frac{1+zs}{z-s}\frac{s^2}{1+s^2}+\frac{s}{1+s^2}
\end{equation}
and
\[\frac{s}{1+s^2}=\frac{s}{1+s^2+t^2}+\frac{s}{1+s^2}\frac{t^2}{1+s^2+t^2}\]
which hold for any $z\in\mathbb{C}\backslash\mathbb{R}$ and $s,t\in\mathbb{R}$, we obtain that
\begin{equation} \label{BFLK1}
\gamma_1+\int_{\mathbb{R}^2}\frac{1+zs}{z-s}\;d\sigma_1(\mathbf{x})=\phi_{\mu^{(1)}}(z)=\mathbf{v}^{(1)}
+\frac{a}{z}+\int_{\mathbb{R}^2}
\left[\frac{zs}{z-s}-\frac{s}{1+s^2+t^2}\right]d\tau(\mathbf{x})
\end{equation}
for some $\mathbf{v}^{(1)}\in\mathbb{R}$. Similarly, one can obtain that
\begin{equation} \label{BFLK2}
\gamma_2+\int_{\mathbb{R}^2}\frac{1+wt}{w-t}\;d\sigma_2(\mathbf{x})=\phi_{\mu^{(2)}}(w)=\mathbf{v}^{(2)}
+\frac{b}{w}+\int_{\mathbb{R}^2}
\left[\frac{wt}{w-t}-\frac{t}{1+s^2+t^2}\right]d\tau(\mathbf{x}).
\end{equation}
Combining the identities (\ref{BFLK1}) and (\ref{BFLK2}) with
\[\widetilde{D}(z,w)=\frac{c}{zw}+\int_{\mathbb{R}^2}\frac{st}{(z-s)(w-t)}\;d\tau(s,t)\] yields the desired expression (\ref{BFLK}). Conversely, a function admitting such an integral form (\ref{BFLK}) with the required properties stated above can be shown to be the bi-free
$\phi$-transform of some bi-freely infinitely divisible measure. These observations lead to the following result.

\begin{thm} \label{bifreeID}
{\normalfont (Bi-free L\'{e}vy-Hin\v{c}in representation)} A probability distribution $\mu$ on $\mathbb{R}^2$ is bi-freely infinitely divisible if and only if its bi-free $\phi$-transform extends analytically to $(\mathbb{C}\backslash\mathbb{R})^2$ and admits the integral representation \eqref{BFLK}, where $\mathbf{v}\in\mathbb{R}^2$, the matrix $\mathbf{A}$ given as in \emph{(\ref{A1})}
is positive semi-definite, and $\tau$ is a positive measure on $\mathbb{R}^2$ with the properties
$\tau(\{\mathbf{0}\})=0$ and $1\wedge\|\mathbf{x}\|^2\in L^1(\tau)$. Moreover, the triplet $(\mathbf{v},\mathbf{A},\tau)$ in \eqref{BFLK} is unique. Conversely, given such a triplet $(\mathbf{v},\mathbf{A},\tau)$ there exists a probability measure $\mu$ for which \eqref{BFLK} holds.
\end{thm}

Theorem \ref{bifreeID} shows that the set $\mathcal{ID}(\bfconvv)$ of bi-freely infinitely divisible distributions is completely parameterized by the triplet $(\mathbf{v},\textbf{A},\tau)$. In the sequel, a probability measure $\mu$ in $\mathcal{ID}(\bfconvv)$ having the representation \eqref{BFLK} will be denoted by $\mu^{(\mathbf{v},\textbf{A},\tau)}_{\bfconvv}$, in which $\tau$ is called the \emph{bi-free L\'{e}vy measure} and $(\mathbf{v},\textbf{A},\tau)$ is called the \emph{bi-free characteristic triplet} of $\mu^{(\mathbf{v},\textbf{A},\tau)}_{\bfconvv}$.
The classical and bi-free L\'evy-Hin\v{c}in representations \eqref{LK1} and \eqref{BFLK} have exactly the same characteristic triplets. As a matter of fact, such a one-to-one correspondence also holds true in general limit theorems, see Theorem \ref{limitthm1}.

\begin{exam} \label{example}
\emph{Let $\mu$ be $\bfconvv$-infinitely divisible and let $(\mathbf{v},\textbf{A},\tau)$ be its bi-free characteristic triplet. \\
(1) Then $\mu$ is called a \emph{bi-free Gaussian distribution} if $\tau=0$. \\
(2) If $\tau$ satisfies
\[\int_{\mathbb{R}^2}\frac{\|\textbf{x}\|}{1+\|\textbf{x}\|^2}\;d\tau(\textbf{x})<\infty,\] then (\ref{BFLK}) is reduced to the form
\[\phi_\mu(z,w)=\frac{\mathbf{u}^{(1)}}{z}+\frac{\mathbf{u}^{(2)}}{w}+\left(\frac{a}{z^2}+\frac{c}{zw}+\frac{b}{w^2}\right)
+\int_{\mathbb{R}^2}\left[\frac{zw}{(z-s)(w-t)}-1\right]d\tau(s,t)\] for some vector $\mathbf{u}\in\mathbb{R}^2$ called the \emph{drift} of $\mu$. \\
(3) Let $\nu$ be in $\mathscr{P}_{\mathbb{R}^2}$ with $\nu(\{\mathbf{0}\})=0$. Then $\mu$ is called a \emph{bi-free compound Poisson distribution} with \emph{rate}
$\lambda>0$ and \emph{jump distribution} $\nu$ if
\begin{equation} \label{compound1}
\mathbf{v}=\int_{\mathbb{R}^2}\frac{\lambda\mathbf{x}}{1+\|\mathbf{x}\|^2}\;d\nu(\mathbf{x}),
\end{equation}
$\mathbf{A}=\mathbf{0}$ and $\tau=\lambda\nu$, and in such a case its bi-free $\phi$-transform is given as
\begin{equation} \label{compound2}
\phi_\mu(z,w)=\lambda\int_{\mathbb{R}^2}\left[\frac{zw}{(z-s)(w-t)}-1\right]d\nu(\mathbf{x}).
\end{equation} If $\nu=\delta_\mathbf{p}$ with $\mathbf{p}\in\mathbb{R}^2\backslash\{\mathbf{0}\}$, then $\mu$ is referred to as a \emph{bi-free Poisson distribution} with rate $\lambda$ and jump distribution $\delta_\mathbf{p}$.}
\end{exam}

\section{Asymptotic behaviors of bi-free convolutions}
This section treats the asymptotic behavior
of the measures
\begin{equation} \label{question}
\mu_n:=\mu_{n1}\bfconv\mu_{n2}\bfconv\cdots\bfconv\mu_{nk_n}\bfconv\delta_{\mathbf{v}_n},
\end{equation}
where $\{k_n\}_{n=1}^\infty$ is a sequence of strictly increasing positive integers, $\{\mu_{nk}\}_{n\geq1,1\leq k\leq k_n}$ is an infinitesimal triangular array in $\mathscr{P}_{\mathbb{R}^2}$, and $\mathbf{v}_n\in\mathbb{R}^2$. In order to cope with the problem, we begin with carrying out the investigation on the asymptotic behavior of the bi-free transforms of $\mu_{nk}$. It turns out that they satisfy certain asymptotic property due to the infinitesimality of $\{\mu_{nk}\}_{n,k}$.

Let $S$ be an unbounded subset of $\mathbb{C}$, and denote by $\mathcal{U}(S^d)$ the collection of
triangular arrays of functions $\{\epsilon_{nk}\}_{n\geq1,1\leq k\leq k_n}$
defined on $S^d$ with the following asymptotic properties: the
functions \[\epsilon_n(z_1,\ldots,z_d)=\max_{1\leq k\leq
k_n}|\epsilon_{nk}(z_1,\ldots,z_d)|,\;\;\;\;\;(z_1,\ldots,z_d)\in S^d,\]
satisfy that $\lim_{n\to\infty}\epsilon_n(z_1,\ldots,z_d)=0$ for any
$(z_1,\ldots,z_d)\in S^d$ and $\epsilon_n(z_1,\ldots,z_d)=o(1)$ uniformly in $n$ as $z_1,\ldots,z_d\to\infty$ with
$(z_1,\ldots,z_d)\in S^d$.

\begin{lem} \label{lemma1} Let $\{\mu_{nk}\}_{n\geq1,1\leq k\leq k_n}\subset\mathscr{P}_{\mathbb{R}^2}$ be infinitesimal. The following statements hold.
\begin{enumerate} [$\qquad(1)$]
\item\label{Existence} {For any $\theta>0$, there exists a number $M>1$ so that
each $\phi_{\mu_{nk}^{(j)}}$ is defined
and satisfies
the relation $\phi_{\mu_{nk}^{(j)}}(\bar{\xi})=\overline{\phi_{\mu_{nk}^{(j)}}(\xi)}$ on $\Gamma:=\Gamma_{\theta,M}$.}
\item\label{Infinitesimal} For $j=1,2,$
\[\left\{\frac{\phi_{\mu_{nk}^{(j)}}(\xi)}{\xi}\right\}_{n\geq1,1\leq k\leq k_n}\in\mathcal{U}(\Gamma)
\;\;\;\;\;\mathrm{and}\;\;\;\;\;\{\phi_{\mu_{nk}}\}_{n\geq1,1\leq
k\leq k_n}\in\mathcal{U}(\Gamma^2).\]
\end{enumerate}
\end{lem}

\begin{pf} For notational convenience, write $G_{nk}=G_{\mu_{nk}}$,
$F_{jnk}=F_{\mu_{nk}^{(j)}}$, $\phi_{jnk} =\phi_{\mu_{nk}^{(j)}}$ and $\phi_{nk}=\phi_{\mu_{nk}}$ for all
$j,k,n$. Let $\theta>0$. Then the existence of the number $M>1$ with the stated properties in \eqref{Existence} is
guaranteed by the infinitesimality of
$\{\mu_{nk}\}_{n,k}$ \cite[Lemma 5]{BerPata00}. The relation $\phi_{jnk}(\bar{z})=\overline{\phi_{jnk}(z)}$ holds for $z\in\Gamma$ because $F_{jnk}(\bar{z})=\overline{F_{jnk}(z)}$ for $z\in \mathbb{C}\backslash\mathbb{R}$.

Recall from \cite[Proposition 2.3]{BerWang} that we can express $\phi_{1nk}$ as
\begin{equation} \label{phiG}
\frac{\phi_{1nk}(z)}{z}=[1+\upsilon_{nk}(z)]\int_\mathbb{R}\frac{s}{z-s}\;d\mu_{nk}^{(1)}(s),
\end{equation}
where $\{\upsilon_{nk}\}_{n\geq1,1\leq k\leq k_n}\in\mathcal{U}(\Gamma)$. Applying the techniques used in the proof of
Proposition \ref{tightG} to the sequence of functions
\[\epsilon_n(z):=\max_{1\leq k\leq k_n}\left|\int_\mathbb{R}\frac{s}{z-s}\;d\mu_{nk}^{(1)}(s)\right|,\] along with the infinitesimality of
$\{\mu_{nk}\}_{n,k}$, yields that $\lim_{n\to\infty}\epsilon_n(z)=0$ for
$z\in\Gamma$ and $\epsilon_n(z)=o(1)$ uniformly in $n$ as $z\to\infty$ with $z\in\Gamma$. Hence the triangular array $\{\phi_{1nk}(z)/z\}_{n,k}$, as well as
$\{\phi_{2nk}(w)/w\}_{n,k}$, is shown to belong to $\mathcal{U}(\Gamma)$.

Finally, in order to show that $\{\phi_{nk}\}_{n\geq1,1\leq k\leq k_n}\in\mathcal{U}(\Gamma^2)$, it suffices to
show that $\{(H_{nk}-1)/H_{nk}\}_{n\geq1,1\leq k\leq k_n}\in\mathcal{U}(\Gamma^2)$,
where $H_{nk}(z,w)=zwG_{nk}(F_{1nk}^{-1}(z),F_{2nk}^{-1}(w))$ for $(z,w)\in\Gamma^2$. This is equivalent to showing that
$\{H_{nk}-1\}_{n,k}\in\mathcal{U}(\Gamma^2)$.
We further see from Proposition \ref{nontangential} that the condition $\{H_{nk}-1\}_{n,k}\in\mathcal{U}(\Gamma^2)$ is the same as \begin{equation} \label{estimatelem1}
\{zwG_{nk}(z,w)-1\}_{n\geq1,1\leq k\leq k_n}\in\mathcal{U}(\Gamma'^2),
\end{equation} where $\Gamma'$ is some truncated cone in $\mathbb{C}$ on which $F_{\mu_{nk}^{(j)}}$ is univalent.
Making use of the techniques in the proof of Proposition \ref{tightG} again one can easily obtain (\ref{estimatelem1}). The proof is complete.
\end{pf} \qed

Let $L>0$ be a fixed number. We will use the following functions and measures to study asymptotic properties of $\mu_n$ defined in (\ref{question}). Let
\begin{equation} \label{vnk}
\mathbf{v}_{nk}=\int_{\{\|\mathbf{x}\|<L\}}\mathbf{x}\;d\mu_{nk}(\mathbf{x})
\end{equation} and define a triangular array
$\{\mathring{\mu}_{nk}\}_{n\geq1,1\leq k\leq k_n}\subset\mathscr{P}_{\mathbb{R}^2}$ as
\begin{equation} \label{barmunk}
\mathring{\mu}_{nk}(B)=\mu_{nk}(B+\mathbf{v}_{nk})
\end{equation} for any Borel set $B\subset\mathbb{R}^2$. This shifted triangular array $\{\mathring{\mu}_{nk}\}$ is also infinitesimal because $\max_{1\leq k \leq k_n}\|\mathbf{v}_{nk}\| \to 0$ as $n\to\infty$. Further define finite positive Borel measures
\begin{equation} \label{nun}
\tau_n=\sum_{k=1}^{k_n}\mathring{\mu}_{nk},
\end{equation}
and functions
\[f_{1nk}(z)=\int_{\mathbb{R}^2}\frac{zs}{z-s}\;d\mathring{\mu}_{nk}(s,t)\;\;\;\;\;
\mathrm{and}\;\;\;\;\;f_{2nk}(w)=\int_{\mathbb{R}^2}\frac{wt}{w-t}\;d\mathring{\mu}_{nk}(s,t)\]
for $z,w\in\mathbb{C}\backslash\mathbb{R}$.

The following result, mostly taken from \cite[Lemma 2.4, Lemma 3.1, Lemma 3.2]{BerWang}, is one of the main ingredients of studying (\ref{question}). For readers' convenience, we provide its proof here.

\begin{lem} \label{fbestimate}
With the same notations $\mu_n$, $\mu_{nk}$ and $\mathbf{v}_n$ in \eqref{question}, $\Gamma$ in {\normalfont Lemma \ref{lemma1}}, and $\tau_n$ and $f_{jnk}$ defined above, the following statements hold.
\begin{enumerate} [$\qquad(1)$]
\item \label{converge-drift} For $j=1,2$ and any fixed $i\ell\in\Gamma$, the sequence
\[\left\{\mathbf{v}_n^{(j)}+\sum_{k=1}^{k_n}[\mathbf{v}_{nk}^{(j)}+f_{jnk}(i\ell)]\right\}_{n=1}^\infty\] converges if and only if the sequence
$\{\phi_{\mu_n^{(j)}}(i\ell)\}_{n=1}^\infty$ converges, in which case they converge to the same value.
\item \label{fjnkestimate}
If
\[V:=\sup_{n\geq1}\int_{\mathbb{R}^2}\frac{\|\mathbf{x}\|^2}{1+\|\mathbf{x}\|^2}
\;d\tau_n(\mathbf{x})<\infty,\]
then there exists an $N\in\mathbb{N}$ such that for $|\ell|\geq1$, $n\geq N$ and $j=1,2$, the inequality
\begin{equation*}
\sum_{k=1}^{k_n}|f_{jnk}(i\ell)|\leq C_LV|\ell|
\end{equation*}
holds for some constant $C_L$ depending only on $L$.
\end{enumerate}
\end{lem}

\begin{pf} We only prove the assertions for $j=1$; the proof for $j=2$ is similar.
Applying the formula (\ref{phiG}) to the triangular array $\{\mathring{\mu}_{nk}^{(1)}\}_{n,k}$, we have
\begin{equation} \label{realation1}
\phi_{\mu_{nk}^{(1)}}(iy)-\mathbf{v}_{nk}^{(1)}=\phi_{\mathring{\mu}_{nk}^{(1)}}(iy)=f_{1nk}(iy)[1+\upsilon_{nk}(iy)]
\end{equation} for $iy\in\Gamma$. Then the desired result in \eqref{converge-drift} follows from \cite[Lemma 2.4]{BerWang} by choosing $z_{nk}=\phi_{\mu_{nk}^{(1)}}(iy)-\mathbf{v}_{nk}^{(1)}$, $w_{nk}=f_{1nk}(iy)$ and $r_n=\mathbf{v}_n^{(1)}+\sum_{k=1}^{k_n}\mathbf{v}_{nk}^{(1)}$.

For assertion \eqref{fjnkestimate}, define
\[b_{1nk}(y)=\mathbf{v}_{nk}^{(1)}\int_{\{\|\mathbf{x}\|\geq L\}}\;d\mu_{nk}(\mathbf{x})+\int_{\{\|\mathbf{x}+\mathbf{v}_{nk}\|\geq L\}}\frac{y^2 s}{y^2+s^2} \;d\mathring{\mu}_{nk}(\mathbf{x}).\] Observe that we have
\[\mathbf{v}_{nk}^{(1)} \int_{\{\|\mathbf{x}\|\geq L\}}\;d\mu_{nk}(\mathbf{x})=\int_{\{\|\mathbf{x}\|<L\}}(s-\mathbf{v}_{nk}^{(1)})\;d\mu_{nk}(\mathbf{x}) =\int_{\{\|\mathbf{x}+\mathbf{v}_{nk}\|<L\}}s\;d\mathring{\mu}_{nk}(\mathbf{x})\]
and
\[f_{1nk}(iy)=\int_{\mathbb{R}^2}\frac{y^2 s}{y^2+s^2}\;d\mathring{\mu}_{nk}(\mathbf{x})
-i\int_{\mathbb{R}^2}\frac{y s^2}{y^2+s^2}\;d\mathring{\mu}_{nk}(\mathbf{x}).\]
This shows that for $|y|>1$,
\begin{equation}\label{bjnkestimate1}
\begin{split}
|\Re f_{1nk}(iy)-b_{1nk}(y)|&=\left|\int_{\{\|\mathbf{x}+\mathbf{v}_{nk}\|<L\}}\frac{s^3}{y^2+s^2}d\mathring{\mu}_{nk}(\mathbf{x})\right|
\leq2L|\Im f_{1nk}(iy)|,
\end{split}
\end{equation}
where we use the fact that $|\mathbf{v}_{nk}^{(1)}|\leq L$ for any $n$ and $k$ in the last inequality.

Now let $N\in\mathbb{N}$ be big enough so that $\sup_{1\leq k \leq k_n}\|\mathbf{v}_{nk}\| \leq L/2$ for all $n\geq N$. Then for $|y|>1$ and $n \geq N$ we have
\begin{equation} \label{bjnkestimate2}
\begin{split}
\sum_{k=1}^{k_n}|b_{1nk}(y)|&\leq\sum_{k=1}^{k_n}\frac{L}{2}\int_{\{\|\mathbf{x}\|\geq L/2\}}\;d\mathring{\mu}_{nk}(\mathbf{x})+
|y|\sum_{k=1}^{k_n}\int_{\{\|\mathbf{x}\|\geq L/2\}} \frac{|y s|}{y^2+s^2}\;d\mathring{\mu}_{nk}(\mathbf{x}) \\
&\leq|y|\sum_{k=1}^{k_n}\int_{\{\|\mathbf{x}\|\geq L/2\}}\frac{1+L}{2} d\mathring{\mu}_{nk}(\mathbf{x}) \\
&=|y|(1+L)(4+L^2)(2L^2)^{-1}\int_{\{\|\mathbf{x}\|\geq L/2\}}\frac{L^2/4}
{1+L^2/4}\;d\tau_n(\mathbf{x}) \\
&\leq|y|(1+L)(4+L^2)(2L^2)^{-1}\int_{\{\|\mathbf{x}\|\geq L/2\}}\frac{\|\mathbf{x}\|^2}
{1+\|\mathbf{x}\|^2}\;d\tau_n(\mathbf{x}).
\end{split}
\end{equation}
Moreover, the estimate
\begin{equation}\label{f1nkestimate}
\sup_{n\geq1}\sum_{k=1}^{k_n}|\Im f_{1nk}(iy)|=
|y|\sup_{n\geq1}\sum_{k=1}^{k_n}\int_{\mathbb{R}^2}\frac{s^2}{y^2+s^2}\;d\mathring{\mu}_{nk}(s,t)\leq V|y|
\end{equation}
holds true for $|y|>1$. Combining \eqref{bjnkestimate1}, \eqref{bjnkestimate2} and \eqref{f1nkestimate} yields assertion \eqref{fjnkestimate}.
\end{pf} \qed

The following result provides an estimation for the bi-free $\phi$-transform of measures $\mu_n$ in (\ref{question}).

\begin{lem} \label{lemma2}
Adopt the notations $\mu_{nk}$, $\mathring{\mu}_{nk}$ and $\Gamma$ in \eqref{question}, \eqref{barmunk} and {\normalfont Lemma \ref{lemma1}}, respectively, and let
\begin{equation} \label{Hnk}
H_{nk}(z,w)=zwG_{\mu_{nk}}\left(F_{\mu_{nk}^{(1)}}^{-1}(z),F_{\mu_{nk}^{(2)}}^{-1}(w)\right),\;\;\;\;\;(z,w)\in\Gamma^2.
\end{equation} Then $H_{nk}-1$ can be expressed as
\[\epsilon_{1nk}\int_{\mathbb{R}^2}\frac{s}{z-s}\;d\mathring{\mu}_{nk}(s,t)+
\epsilon_{2nk}\int_{\mathbb{R}^2}\frac{t}{w-t}\;d\mathring{\mu}_{nk}(s,t)+(1+\epsilon_{nk})
\int_{\mathbb{R}^2}\frac{st}{(z-s)(w-t)}\;d\mathring{\mu}_{nk}(s,t),\]
where $\{\epsilon_{1nk}(z,w)\}_{n,k}$, $\{\epsilon_{2nk}(z,w)\}_{n,k}$ and $\{\epsilon_{nk}(z,w)\}_{n,k}$ are
triangular arrays of functions in $\mathcal{U}(\Gamma^2)$. Consequently, $\{H_{nk}-1\}\in\mathcal{U}(\Gamma^2)$.
\end{lem}

\begin{pf} For notational convenience, write $\mathring{G}_{nk}=G_{\mathring{\mu}_{nk}}$, $\mathring{F}_{jnk}=F_{\mathring{\mu}_{nk}^{(j)}}$ and
$\mathring{\phi}_{jnk}=\phi_{\mathring{\mu}_{nk}^{(j)}}$ for any $j$, $n$ and $k$. Then a simple argument of change of variables shows that
\[H_{nk}(z,w)=zw\mathring{G}_{nk}\big(\mathring{F}_{1nk}^{-1}(z),\mathring{F}_{2nk}^{-1}(w)\big),\;\;\;\;\;(z,w)\in\Gamma^2.\]
To conclude the proof, we first analyze the function
\[H_{nk}(z,w)-1=\int_{\mathbb{R}^2}\left[\frac{z}{\mathring{\phi}_{1nk}(z)+z-s}\frac{w}{\mathring{\phi}_{2nk}(w)+w-t}-1\right]
d\mathring{\mu}_{nk}(s,t).\]
Using the identity
\[\frac{\xi}{\mathring{\phi}_{jnk}(\xi)+\xi-r}=\frac{\xi}{\xi-r}\left[1-\frac{\mathring{\phi}_{jnk}(\xi)}{\mathring{F}_{jnk}^{-1}(\xi)-r}\right],\]
the function $H_{nk}-1$ can be rewritten as the sum of functions $I_{1nk}$,
$I_{2nk}$ and $I_{3nk}$, where
\[I_{1nk}(z,w)=\int_{\mathbb{R}^2}\left[\frac{zw}{(z-s)(w-t)}-1\right]d\mathring{\mu}_{nk}(s,t),\]
\[I_{2nk}(z,w)=-\int_{\mathbb{R}^2}\frac{zw}{(z-s)(w-t)}\left[\frac{\mathring{\phi}_{1nk}(z)}{\mathring{F}_{1nk}^{-1}(z)-s}+
\frac{\mathring{\phi}_{2nk}(w)}{\mathring{F}_{2nk}^{-1}(w)-t}\right]d\mathring{\mu}_{nk}(s,t),\] and
\[I_{3nk}(z,w)=\int_{\mathbb{R}^2}\frac{zw}{(z-s)(w-t)}\frac{\mathring{\phi}_{1nk}(z)}{\mathring{F}_{1nk}^{-1}(z)-s}
\frac{\mathring{\phi}_{2nk}(w)}{\mathring{F}_{2nk}^{-1}(w)-t}\;d\mathring{\mu}_{nk}(s,t).\]
For any $z\in\Gamma$ with $|z|$ large enough, $s\in\mathbb{R}$, and $n,k$, Lemma \ref{lemma1} shows that
\[\left|\frac{\mathring{F}_{1nk}^{-1}(z)-s}{z}\right|=\left|\frac{z-s}{z}+\frac{\mathring{\phi}_{1nk}(z)}{z}\right|\geq
\frac{|\Im
z|}{|z|}-\left|\frac{\mathring{\phi}_{1nk}(z)}{z}\right|\geq\frac{1}{\sqrt{1+\theta^2}}-\left|\frac{\mathring{\phi}_{1nk}(z)}{z}\right|\geq
c_\theta,\] where $c_\theta=1/(2\sqrt{1+\theta^2})$. This implies
that
\begin{equation} \label{lemma2eq1}
\frac{1}{\mathring{F}_{1nk}^{-1}(z)-s}=\frac{1}{z-s}\left[1-\frac{\mathring{\phi}_{1nk}(z)}{\mathring{F}_{1nk}^{-1}(z)-s}\right]=\frac{1}{z-s}
[1+\delta_{1nk}(z,s)],
\end{equation}
where
\[\max_{1\leq k\leq k_n}|\delta_{1nk}(z,s)|=\max_{1\leq k\leq k_n}\left|\frac{\mathring{\phi}_{1nk}(z)}{\mathring{F}_{1nk}^{-1}(z)-s}\right|
\leq c_\theta^{-1}\max_{1\leq k\leq
k_n}\left|\frac{\mathring{\phi}_{1nk}(z)}{z}\right|:=\delta_{1n}(z).\]
We further obtain from the estimate (\ref{lemma2eq1}) that
\begin{equation} \label{estimateI2nk}
\begin{split}
&\max_{1\leq k\leq k_n}\left|\int_{\mathbb{R}^2}\frac{z^2w}{(z-s)(w-t)}\frac{1}{\mathring{F}_{1nk}^{-1}(z)-s}\;d\mathring{\mu}_{nk}(s,t)-1\right| \\
=&\max_{1\leq k\leq k_n}\left|\int_{\mathbb{R}^2}\frac{z^2w}{(z-s)^2(w-t)}[1+\delta_{1nk}(z,s)]\;d\mathring{\mu}_{nk}(s,t)-1\right| \\
=&\max_{1\leq k\leq k_n}\left|\int_{\mathbb{R}^2}\left\{\left[\frac{z^2w}{(z-s)^2(w-t)}-1\right][1+\delta_{1nk}(z,s)]+\delta_{1nk}(z,s)
\right\}d\mathring{\mu}_{nk}(s,t)\right| \\
\leq&\delta_{1n}(z)+[1+\delta_{1n}(z)]M_n(z,w),
\end{split}
\end{equation} where
\[M_n(z,w)=\max_{1\leq k\leq k_n}\int_{\mathbb{R}^2}\left|\frac{z^2w}{(z-s)^2(w-t)}-1\right|d\mathring{\mu}_{nk}(s,t).\]
Note that the function
$\delta_{1n}(z)+[1+\delta_{1n}(z)]M_n(z,w)=o(1)$ as $n\to\infty$ for $(z,w)\in\Gamma^2$ and uniformly in $n$ as $z,w\to\infty$ with $(z,w)\in\Gamma^2$. Indeed, this can be easily
obtained by using Lemma \ref{lemma1} and applying the techniques
employed in the proof of Proposition \ref{tightG} to the identity
\[\frac{z^2w}{(z-s)^2(w-t)}-1=\frac{z}{z-s}\left[\frac{zw}{(z-s)(w-t)}-1\right]+\frac{s}{z-s}.\]
Now applying the formula \eqref{phiG} to the triangular array $\{\mathring{\mu}_{nk}\}_{n\geq1,1\leq k\leq k_n}$ and using \eqref{estimateI2nk} give
\begin{align*}
\int_{\mathbb{R}^2}\frac{zw}{(z-s)(w-t)}\frac{\mathring{\phi}_{1nk}(z)}{\mathring{F}_{1nk}^{-1}(z)-s}\;d\mathring{\mu}_{nk}(s,t)&=
\frac{\mathring{\phi}_{1nk}(z)}{z}\int_{\mathbb{R}^2}\frac{z^2w}{(z-s)(w-t)}\frac{d\mathring{\mu}_{nk}(s,t)}{\mathring{F}_{1nk}^{-1}(z)-s} \\
&=[1+v_{1nk}(z,w)]\int_{\mathbb{R}^2}\frac{s}{z-s}\;d\mathring{\mu}_{nk}(s,t),
\end{align*} where $\{v_{1nk}\}_{n,k}$ is a triangular array in $\mathcal{U}(\Gamma^2)$. Similarly, the identity
\[\int_{\mathbb{R}^2}\frac{zw}{(z-s)(w-t)}\frac{\mathring{\phi}_{2nk}(w)}{\mathring{F}_{2nk}^{-1}(w)-t}\;d\mathring{\mu}_{nk}(s,t)
=[1+v_{2nk}(z,w)]\int_{\mathbb{R}^2}\frac{t}{w-t}\;d\mathring{\mu}_{nk}(s,t)\]
is valid for some triangular array
$\{v_{2nk}\}_{n,k}\in\mathcal{U}(\Gamma^2)$. By similar arguments, one can also show that
\[I_{3nk}=v_{3nk}\int_{\mathbb{R}^2}\frac{s}{z-s}\;d\mathring{\mu}_{nk}(s,t)+
v_{4nk}\int_{\mathbb{R}^2}\frac{t}{w-t}\;d\mathring{\mu}_{nk}(s,t)
+v_{5nk}\int_{\mathbb{R}^2}\frac{st\;d\mathring{\mu}_{nk}(s,t)}{(z-s)(w-t)}\] for some triangular arrays $\{v_{3nk}\}_{n,k}$, $\{v_{4nk}\}_{n,k}$ and $\{v_{5nk}\}_{n,k}$ in $\mathcal{U}(\Gamma^2)$. Finally, let $\epsilon_{1nk}=v_{3nk}-v_{1nk}$, $\epsilon_{2nk}=v_{4nk}-v_{2nk}$, and $\epsilon_{nk}=v_{5nk}$. Then we conclude the proof of the first assertion by
using the integral representations of $I_{2nk}$ and $I_{3nk}$ provided above and the identity
\[I_{1nk}(z,w)=\int_{\mathbb{R}^2}\left[\frac{s}{z-s}+\frac{t}{w-t}+\frac{st}{(z-s)(w-t)}\right]d\mathring{\mu}_{nk}(s,t).
\] The fact that $\{H_{nk}-1\}\in\mathcal{U}(\Gamma^2)$ can be proved by the infinitesimality of $\{\mathring{\mu}_{nk}\}$ and the techniques in Proposition \ref{tightG}. This finishes the proof.
\end{pf} \qed

\begin{lem} \label{lemma3}
Suppose that the marginal laws of $\mu_n$ in \eqref{question} converge weakly. With the notations $\mathring{\mu}_{nk}$, $\tau_n$, $\Gamma$ as before and $H_{nk}$ as in \emph{(\ref{Hnk})}, the following statements hold.
\begin{enumerate} [$\qquad(1)$]
\item\label{conv-sigma} {The positive planar measures $\{\sigma_{1n}\}_{n=1}^\infty$ and $\{\sigma_{2n}\}_{n=1}^\infty$ defined as
\[\sigma_{1n}=\frac{s^2}{1+s^2}\;\tau_n\;\;\;\;\;\mathrm{and}\;\;\;\;\;\sigma_{2n}=\frac{t^2}{1+t^2}\;\tau_n\]
are uniformly bounded and tight.}
\item\label{convH} {For $(iy,iv)\in\Gamma^2$, the limit
\begin{equation} \label{limit1}
\lim_{n\to\infty}\sum_{k=1}^{k_n}\frac{H_{nk}(iy,iv)-1}{H_{nk}(iy,iv)}
\end{equation} exists if and only if the limit
\begin{equation} \label{limit2}
\lim_{n\to\infty}\sum_{k=1}^{k_n}\int_{\mathbb{R}^2}\frac{st}{(iy-s)(iv-t)}\;d\mathring{\mu}_{nk}(s,t)
\end{equation} exists, in which case these two limits are equal.}
\item\label{estimateH} {The function
\[\sum_{k=1}^{k_n}\frac{H_{nk}(iy,iv)-1}{H_{nk}(iy,iv)}=o(1)\] uniformly in $n$ as $|y|,|v|\to\infty$ if and only if
\begin{equation} \label{o(1)}
\sum_{k=1}^{k_n}\int_{\mathbb{R}^2}\frac{st}{(iy-s)(iv-t)}\;d\mathring{\mu}_{nk}(s,t)=o(1)
\end{equation}
uniformly in $n$ as $|y|,|v|\to\infty$.}
\end{enumerate}
\end{lem}

\begin{pf} For any $\epsilon>0$, choose a large positive number $y_0>1$ so that $|\phi_{\mu_n^{(1)}}(i y_0)|<\epsilon y_0$ for all $n$ by Proposition \ref{weakconvphi}. Then we deduce from Lemma \ref{fbestimate}\eqref{converge-drift} and the identity
\begin{align*}
\int_{\mathbb{R}^2}\frac{s^2}{y_0^2+s^2}\;d\tau_n(s,t)=- \frac{1}{y_0}\Im\left(\mathbf{v}_n^{(1)}+\sum_{k=1}^{k_n}[\mathbf{v}_{nk}^{(1)}+f_{1nk}(iy_0)]\right)
\end{align*}
the existence of a large number $N\in\mathbb{N}$ so that
\[\int_{\mathbb{R}^2}\frac{s^2}{y_0^2+s^2}\;d\tau_n(s,t)<2\epsilon,\;\;\;\;\;n\geq N.\]
This, along with the inequalities
\[\frac{s^2}{1+s^2}\leq\frac{y_0^2s^2}{y_0^2+s^2}\;\;\;\;\;\mathrm{and}\;\;\;\;\;\frac{s^2}{1+s^2}\leq\frac{2s^2}{y_0^2+s^2}\] which hold true for $s\in\mathbb{R}$ and $|s|\geq y_0$, respectively, yields the uniform boundedness and tightness of $\{\sigma_{1n}\}$. Similarly, $\{\sigma_{2n}\}$ is uniformly bounded and tight. This proves (1).

To prove \eqref{convH}, we first argue that the limit in (\ref{limit2}) exists if and only if so does the limit
\begin{equation} \label{limit3}
\lim_{n\to\infty}\sum_{k=1}^{k_n}[H_{nk}(iy,iv)-1],
\end{equation} and show that they are equal. We shall use the integral representation of $H_{nk}-1$ given in Lemma \ref{lemma2} to accomplish this goal. Observe first that the quantity $\V$ defined in Lemma \ref{fbestimate}\eqref{fjnkestimate} is finite by the established result \eqref{conv-sigma}. Hence for $(iv,iy)\in\Gamma^2$ and all large $n$, we have
\begin{align*}
\left|\sum_{k=1}^{k_n}\epsilon_{1nk}(iy,iv)\int_{\mathbb{R}^2}\frac{s}{iy-s}\;d\mathring{\mu}_{nk}(s,t)\right|
&\leq\left[\max_{1\leq k\leq k_n}|\epsilon_{1nk}(iy,iv)|\right]\frac{1}{|y|}\sum_{k=1}^{k_n}|f_{1nk}(iy)| \\
&\leq C_L \V \max_{1\leq k\leq k_n}|\epsilon_{1nk}(iy,iv)|,
\end{align*} which yields that
\[\lim_{n\to\infty}\sum_{k=1}^{k_n}\epsilon_{1nk}(iy,iv)\int_{\mathbb{R}^2}\frac{s}{iy-s}\;d\mathring{\mu}_{nk}(s,t)=0,\] as well as
\[\lim_{n\to\infty}\sum_{k=1}^{k_n}\epsilon_{2nk}(iy,iv)\int_{\mathbb{R}^2}\frac{t}{iv-t}\;d\mathring{\mu}_{nk}(s,t)=0\]
by similar arguments. Next notice that the inequality
\[\left|\frac{st}{(iy-s)(iv-t)}\right|\leq\frac{2s^2}{1+s^2}+\frac{2 t^2}{1+t^2},\;\;\;\;\;|y|,|v|\geq1,\;s,t\in\mathbb{R},\] along with the established result (1), allows us to obtain that
\[\lim_{n\to\infty}\sum_{k=1}^{k_n}\epsilon_{nk}(iy,iv)
\int_{\mathbb{R}^2}\frac{st}{(iy-s)(iv-t)}\;d\mathring{\mu}_{nk}(s,t)=0\]
for any point $(iy,iv)\in\Gamma^2$. Hence we have proved that the pointwise convergence of (\ref{limit2}) is equivalent to that of (\ref{limit3}), and both limits are the same.
These discussions also indicate that $\sum_{k=1}^{k_n}|H_{nk}(iy,iv)-1|$ is uniformly bounded in $n$ for fixed $(iy,iv)\in\Gamma^2$. Since $\{H_{nk}-1\}_{nk}\in\mathcal{U}(\Gamma^2)$, one can see that (\ref{limit1}) converges pointwise if and only if so does (\ref{limit3}), and they have the same limit. This yields assertion \eqref{convH}. The preceding discussions with a little effort yield assertion \eqref{estimateH}.
\end{pf} \qed

Before stating the main theorem of this section, let us introduce the following conditions, which play an important role in the asymptotic problem under investigation.

\begin{cond} \label{cond1}
Let $\{\tau_n\}_{n=1}^\infty$ be a sequence of finite positive Borel measures on $\mathbb{R}^2$.
\begin{enumerate}[\quad(I)]
\item\label{I} The sequences of measures
\begin{equation*} \label{sigmajn}
\sigma_{1n}=\frac{s^2}{1+s^2}\;\tau_n\;\;\;\;\;
\mathrm{and}\;\;\;\;\;
\sigma_{2n}=\frac{t^2}{1+t^2}\;\tau_n
\end{equation*} converge weakly to finite positive Borel measures $\sigma_1$ and $\sigma_2$ on $\mathbb{R}^2$, respectively.
\item\label{II} The limit
\[\gamma:=\lim_{n\to\infty}\int_{\mathbb{R}^2}\frac{st}{(1+s^2)(1+t^2)}\;d\tau_n(s,t)\]
exists in $\mathbb{R}$.
\end{enumerate}
\end{cond}

\begin{thm} \label{nonidentical}
Let $\{\mathbf{v}_n\}_{n=1}^\infty$ be a sequence of vectors in $\mathbb{R}^2$, $\{k_n\}_{n=1}^\infty$ a sequence of strictly increasing positive integers, and let
$\{\mu_{nk}\}_{n\geq1,1\leq k\leq k_n}$ be an infinitesimal triangular array in $\mathscr{P}_{\mathbb{R}^2}$. Following the notations in \eqref{vnk}, \eqref{barmunk} and \eqref{nun}, the following statements are equivalent.
\begin{enumerate} [$\qquad(1)$]

\item\label{Statement1}
The sequence
\[\mu_n:=\mu_{n1}\bfconv\mu_{n2}\bfconv\cdots\bfconv\mu_{nk_n}\bfconv\delta_{\mathbf{v}_n}\]
converges weakly to some planar probability measure $\mu_{\bfconvv}$.

\item\label{Statement2}
{\normalfont Condition \ref{cond1}}\eqref{I} and \eqref{II} are satisfied, and the sequence
\begin{equation} \label{vector}
\mathbf{v}_n+\sum_{k=1}^{k_n}\left[\mathbf{v}_{nk}+\int_{\mathbb{R}^2}\frac{\mathbf{x}}{1+\|\mathbf{x}\|^2}\;d\mathring{\mu}_{nk}(\mathbf{x})\right]
\end{equation}
of vectors in $\mathbb{R}^2$ converges to some vector $\mathbf{v}$.
\end{enumerate}
\end{thm}

\begin{pf} We take the set $\Gamma$ given in Lemma \ref{lemma1}. Suppose that assertion \eqref{Statement1} holds. By Lemma \ref{lemma3}\eqref{conv-sigma}, let $\sigma_{1j_n}\Rightarrow\sigma_1$ and
$\sigma_{2j_n}\Rightarrow\sigma_2$ for some subsequences $\{\sigma_{1j_n}\},\{\sigma_{2j_n}\}$ and for some finite positive Borel measures $\sigma_1,\sigma_2$. Observe next that the limit in (\ref{limit2}) exists. Denote this limit by $\widetilde{K}(iy,iv)$. Then using the
decomposition (\ref{decom}) we see that that the limit
\[\gamma':=\lim_{n\to\infty}\int_{\mathbb{R}^2}\frac{st}{(1+s^2)(1+t^2)}\;d\tau_{j_n}(s,t)\] must exist
and $\widetilde{K}$ has an analytic extension (still denoted by $\widetilde{K}$) to $(\mathbb{C}\backslash\mathbb{R})^2$. More precisely, the analytic extension $K(z,w):=zw\widetilde{K}(z,w)$ can be expressed as the sum of integrals:
\[K(z,w)=\gamma'+\int_{\mathbb{R}^2}\frac{1+zs}{z-s}\left[\frac{t}{1+t^2}
+\frac{1+wt}{w-t}\frac{t^2}{1+t^2}\right]d\sigma_1(s,t)+
\int_{\mathbb{R}^2}\frac{1+wt}{w-t}\frac{s}{1+s^2}\;d\sigma_2(s,t).\]
A simple calculation then shows that for any $z=x+iy$ and $w=u+iv$ in $\mathbb{C}^+$,
\[-\frac{1}{2}\Re[K(z,w)-K(\overline{z},w)]=\int_{\mathbb{R}^2}\frac{yv}{[(s-x)^2+y^2][(t-u)^2+v^2]}(1+s^2)t^2\;d\sigma_1(s,t).\]
This identity, of course, is also valid for any other weak-limit point $\sigma_1'$ of $\{\sigma_{1n}\}$.
Combining this result with the Stieltjes inversion formula for two variables shows that
\begin{equation} \label{unique1}
\frac{t^2}{1+t^2}\;\sigma_1=\frac{t^2}{1+t^2}\;\sigma_1'.
\end{equation}

On the other hand, according to Lemma \ref{fbestimate}\eqref{converge-drift} the sequence
\begin{equation} \label{gamma1n}
\gamma_{1n}:=\mathbf{v}_n^{(1)}+\sum_{k=1}^{k_n}\left[\mathbf{v}_{nk}^{(1)}
+\int_{\mathbb{R}^2}\frac{s}{1+s^2}\;d\mathring{\mu}_{nk}(s,t)\right],
\end{equation} converges to some $\gamma_1$ and
\[\phi_{\mu_{\bfconvv}^{(1)}}(z)=\gamma_1+\int_{\mathbb{R}^2}\frac{1+zs}{z-s}\;d\sigma_1(s,t),\;\;\;\;\;
z\in\mathbb{C}\backslash\mathbb{R}.\]
Hence $\mu_{\bfconvv}^{(1)}$ is $\boxplus$-infinitely divisible and $(\gamma_1,\sigma_1^{(1)})$, as well as $(\gamma_1,\sigma_1'^{(1)})$, is the free generating pair for $\mu_{\bfconvv}^{(1)}$ \cite{BerPata99}.
Therefore we obtain that $\sigma_1^{(1)}=\sigma_1'^{(1)}$. This with (\ref{unique1}) shows that $\sigma_1=\sigma_1'$ by \cite[Lemma 3.10]{HW}. We conclude that $\sigma_1$ is the unique weak-limit point of $\{\sigma_{1n}\}_{n=1}^\infty$. Similarly, $\{\sigma_{2n}\}$ has only one weak-limit point. Hence Condition \ref{cond1}\eqref{I} and \ref{cond1}\eqref{II} are satisfied. Moreover, the identity
\[\frac{s}{1+s^2}-\frac{s}{1+s^2+t^2}=\frac{st^2}{(1+s^2+t^2)(1+s^2)}\]
shows that the vector defined in (\ref{vector}) converges. Hence the proof of $\eqref{Statement1}\Rightarrow\eqref{Statement2}$ is complete.

Conversely, suppose that assertion \eqref{Statement2} holds. Then $\mathbf{v}_n^{(j)}+\sum_{k=1}^{k_n}[\mathbf{v}_{nk}^{(j)}+f_{jnk}(i\ell)]$ converges as $n\to\infty$ for any $i\ell\in\Gamma$ and $j=1,2$, and hence so does $\phi_{\mu_n^{(j)}}(i\ell)$ by Lemma \ref{fbestimate}\eqref{converge-drift}.
Employing the identity (\ref{decom}) gives that the
limit (\ref{limit2}) must exist, from which we see that $\phi_{\mu_n}$ converges
pointwise on $\Gamma^2$ by Lemma \ref{lemma3}\eqref{convH}. To finish the proof of $\eqref{Statement2}\Rightarrow\eqref{Statement1}$,
it remains to show that $\phi_{\mu_n^{(1)}}(iy)=o(|y|)$, $\phi_{\mu_n^{(2)}}(iv)=o(|v|)$ and
$\sum_{k=1}^{k_n}[H_{nk}(iy,iv)-1]/H_{nk}(iy,iv)=o(1)$
uniformly in $n$ as $|y|,|v|\to\infty$ by Proposition
\ref{weakconvphi}.
First of all, the identities (\ref{decom}) and (\ref{realation1}) show that
\begin{equation} \label{phiunif}
\phi_{\mu_n^{(1)}}(iy)=\gamma_{1n}+\int_{\mathbb{R}^2}\frac{1+iys}{iy-s}\;d\sigma_{1n}(s,t)+
\sum_{k=1}^{k_n}\upsilon_{nk}(iy)f_{1nk}(iy),
\end{equation}
where $\gamma_{1n}$ is defined as in (\ref{gamma1n}) and $\upsilon_{nk}\in\mathcal{U}(\Gamma)$.
Since $\{\gamma_{1n}\}$ and $\{y^{-1}\sum_{k=1}^{k_n}|f_{1nk}(iy)|:iy\in\Gamma\}$ are bounded for all large $n$ by the hypotheses in \eqref{Statement2} and Lemma \ref{fbestimate}\eqref{fjnkestimate}, it
suffices to show that the second term in (\ref{phiunif}) equals $o(|y|)$ uniformly in $n$ as $|y|\to\infty$.
Notice that for any $r>0$ and $|y|\geq1$, we have
\begin{align*}
\frac{1}{|y|}\int_{\mathbb{R}^2}\left|\frac{1+iys}{iy-s}\right|d\sigma_{1n}(s,t)&\leq\frac{1}{|y|}
\int_{\{\|\mathbf{x}\|\leq r\}}\frac{1+|sy|}{\sqrt{y^2+s^2}}\;d\sigma_{1n}(s,t)+\sigma_{1n}(\{\|\mathbf{x}\|>r\}) \\
&\leq\frac{1}{|y|^2}(1+r|y|)\sup_n\sigma_{1n}(\mathbb{R}^2)+\sigma_{1n}(\{\|\mathbf{x}\|>r\}),
\end{align*} which yields the desired result by the uniform boundedness and tightness of $\{\sigma_{1n}\}$. Similarly, $\phi_{\mu_n^{(2)}}(iv)=o(|v|)$ uniformly in $n$ as $|v|\to\infty$. The last desired result is equivalent to
the uniform convergence of (\ref{o(1)}) in $n$ as $|y|,|v|\to\infty$ by Lemma \ref{lemma3}\eqref{estimateH}. The latter uniform convergence can be proved by using Condition \ref{cond1}(II) and applying the techniques shown above to the integral in (\ref{o(1)}), which is rewritten as
\[-\frac{1}{yv}\sum_{k=1}^{k_n}\int_{\mathbb{R}^2}\left[\frac{1+iys}{iy-s}\frac{s^2}{1+s^2}
+\frac{s}{1+s^2}
\right]\left[\frac{1+ivt}{iv-t}\frac{t^2}{1+t^2}+\frac{t}{1+t^2}
\right]d\mathring{\mu}_{nk}(s,t).\] This finishes the proof of $(2)\Rightarrow(1)$.
\end{pf} \qed

Let $\mathbf{u}_n$ be the vector defined in (\ref{vector}). From the proof of Theorem \ref{nonidentical}, one can see that for large $|y|$, $\phi_{\mu_n^{(1)}}(iy)$
can also be expressed as
\[\phi_{\mu_n^{(1)}}(iy)=\mathbf{u}_n^{(1)}+\sum_{k=1}^{k_n}\left[\big(1+\upsilon_{nk}(iy)\big)f_{1nk}(iy)-\int_{\mathbb{R}^2}
\frac{s}{1+s^2+t^2}\;d\mathring{\mu}_{nk}(s,t)\right].\] A similar expression also holds for $\phi_{\mu_n^{(2)}}(iv)$ when $|v|$ is large.
Since
\[\phi_{\mu_{\bfconvv}}(iy,iv)=\lim_{n\to\infty}\left[
\frac{\phi_{\mu_n^{(1)}}(iy)}{iy}+\frac{\phi_{\mu_n^{(2)}}(iv)}{iv}+\sum_{k=1}^{k_n}\int_{\mathbb{R}^2}\frac{st}{(iy-s)(iv-t)}
\;d\mathring{\mu}_{nk}(s,t)\right],\] a simple calculation shows that
\begin{equation} \label{BFLKyv}
\phi_{\mu_{\bfconvv}}(iy,iv)=\frac{\mathbf{v}^{(1)}}{iy}+\frac{\mathbf{v}^{(2)}}{iv}+\mathcal{S}(iy,iv)
\end{equation} for $(iy,iv)\in\Gamma^2$, where
\[\mathcal{S}(iy,iv)=\lim_{n\to\infty}
\sum_{k=1}^{k_n}\int_{\mathbb{R}^2}\left[\frac{(iy)(iv)}{(iy-s)(iv-t)}
-1-\frac{(iy)^{-1}s+(iv)^{-1}t}{1+s^2+t^2}\right]d\mathring{\mu}_{nk}(s,t).\] In the next section we shall use (\ref{BFLKyv}) to show that the function $\phi_{\mu_{\bfconvv}}$ extends analytically to $(\mathbb{C}\backslash\mathbb{R})^2$ and the analytic extension admits an integral representation of the from (\ref{BFLK}). As a consequence of Theorem \ref{bifreeID}, $\mu_{\bfconvv}$ is bi-freely infinitely divisible.

\section{Transfer principle for limit theorems and bijection between $\mathcal{ID}(*)$ and $\mathcal{ID}(\bfconvv)$}
This section is mainly devoted to studying the relation between the sets $\mathcal{ID}(*)$ and $\mathcal{ID}(\bfconvv)$ in terms of classical and bi-free limit theorems. We first introduce another type of convergence on the set of positive Borel measures on $\mathbb{R}^2$.

\begin{pdef} {\normalfont Denote by $\mathcal{M}_{\mathbb{R}^2}^\mathbf{0}$ the set of positive Borel (not necessarily finite) measures $\tau$ on $\mathbb{R}^2$ for which $\tau(B)<\infty$ for any Borel set $B\subset\mathbb{R}^2$ bounded away from zero, i.e. $\inf_{x\in B}\|x\|>0$. For a measure $\tau$ and a sequence of measures $\{\tau_n\}_{n=1}^\infty$ in $\mathcal{M}_{\mathbb{R}^2}^\mathbf{0}$, the convergent situation that $\tau_n(B)\to\tau(B)$ for any Borel set $B$ which is bounded away from zero and satisfies $\tau(\partial B)=0$ is denoted by $\tau_n\Rightarrow_\mathbf{0}\tau$.}
\end{pdef}

We remark here that any finite positive Borel measure on $\mathbb{R}^2$ belongs to $\mathcal{M}_{\mathbb{R}^2}^\mathbf{0}$.
Also note that the limiting measure $\tau$ in the convergence $\tau_n\Rightarrow_\mathbf{0}\tau$ is not necessarily unique since an arbitrary mass at $\mathbf{0}$ can be added to it. Portmanteau theorem for measures in $\mathcal{M}_{\mathbb{R}^2}^\mathbf{0}$ is stated below  (see \cite{Portmanteau06}).

\begin{prop} \label{Portmanteau}
Given $\{\tau_n\}_{n=1}^\infty$ and $\tau$ in $\mathcal{M}_{\mathbb{R}^2}^\mathbf{0}$, the following are equivalent:
\begin{enumerate} [$\qquad(1)$]
\item {$\tau_n\Rightarrow_\mathbf{0}\tau$;}
\item {for any bounded and continuous function $f$ on $\mathbb{R}^2$ with support bounded away from zero,
\[\lim_{n\to\infty}\int_{\mathbb{R}^2}f\;d\tau_n=\int_{\mathbb{R}^2}f\;d\tau;\]}
\item {for any bounded and continuous function $f$ on $\mathbb{R}^2$ and for any
Borel set $B\subset\mathbb{R}^2$ which is bounded away from zero and satisfies $\tau(\partial B)=0$,
\[\lim_{n\to\infty}\int_Bf\;d\tau_n=\int_Bf\;d\tau;\]}
\item\label{Port4} for every closed subset $C$ and open subset $O$ of $\mathbb{R}^2$
that are bounded away from zero,
\[\limsup_{n\to\infty}\tau_n(C)\leq\tau(C)\text{\qquad and\qquad}\liminf_{n\to\infty}\tau_n(O)\geq\tau(O).\]
\end{enumerate}
\end{prop}

We next introduce two conditions that are used in characterizing the classical limit theorem in multidimensional spaces (see \cite{limitthm}).

\begin{cond} \label{cond2}
Let $\{\tau_n\}_{n=1}^\infty$ be a sequence of measures in $\mathcal{M}_{\mathbb{R}^2}^\mathbf{0}$.
\begin{enumerate}[\quad(I)]
\setcounter{enumi}{2}
\item\label{III} $\tau_n\Rightarrow_\mathbf{0}\tau$ for some measure
$\tau\in\mathcal{M}_{\mathbb{R}^2}^\mathbf{0}$ with $\tau(\{\mathbf{0}\})=0$;
\item\label{IV} for every vector $\mathbf{u}\in\mathbb{R}^2$, the limits
\begin{equation} \label{limsupinf}
\lim_{\epsilon\to0^+}\limsup_{n\to\infty}\int_{\|\mathbf{x}\|<\epsilon}
\langle\mathbf{u},\mathbf{x}\rangle^2d\tau_n(\mathbf{x})\quad\mathrm{and}\quad
\lim_{\epsilon\to0^+}\liminf_{n\to\infty}\int_{\|\mathbf{x}\|<\epsilon}
\langle\mathbf{u},\mathbf{x}\rangle^2d\tau_n(\mathbf{x})
\end{equation} exist (as finite numbers), and they are equal.
\end{enumerate}
\end{cond}

In the following, we show the equivalence between Condition \ref{cond1} and Condition \ref{cond2}, which will play an important role in clarifying the relation between $\mathcal{ID}(*)$ and $\mathcal{ID}(\bfconvv)$.

\begin{lem} \label{mainlem}
Let $\{\tau_n\}_{n=1}^\infty$ be a sequence of finite positive Borel measures on $\mathbb{R}^2$. Then \eqref{I} and \eqref{II} in {\normalfont Condition \ref{cond1}} hold if and only if \eqref{III} and \eqref{IV} in {\normalfont Condition \ref{cond2}} are satisfied, in which case
\begin{equation} \label{quantity}
c=\gamma-\int_{\mathbb{R}^2}\frac{st}{(1+s^2)(1+t^2)}\;d\tau(s,t),
\end{equation} is a finite number, the matrix
\begin{equation} \label{matrix}
\mathbf{A}=\left(
\begin{array}{cc}
\sigma_1(\{\mathbf{0}\}) & c \\
c & \sigma_2(\{\mathbf{0}\}) \\
\end{array}
\right)
\end{equation}
is positive semi-definite and the limits in \eqref{limsupinf} define a non-negative quadratic form $\langle\mathbf{A}\mathbf{u},\mathbf{u}\rangle$.
\end{lem}

\begin{pf} Suppose that \eqref{I} and \eqref{II} in Condition \ref{cond1} are satisfied. With $\sigma_1$ and $\sigma_2$ at hand, one can define the positive measure $\tau$ as in (\ref{taudef}) with $\tau(\{\mathbf{0}\})$. Then the relation
\[\frac{t^2}{1+t^2}\;\sigma_1=\frac{s^2}{1+s^2}\;\sigma_2,\] which is obtained from the definition of $\sigma_{1n}$ and $\sigma_{2n}$, ensures that $\tau$ is well defined. It is also easy to verify that $\tau(\{\|\mathbf{x}\|\geq\epsilon\})<\infty$ for any $\epsilon>0$, whence $\tau\in\mathcal{M}_{\mathbb{R}^2}^\mathbf{0}$. Now we claim $\tau_n\Rightarrow_\mathbf{0}\tau$.
Pick any bounded and continuous function $f$ on $\mathbb{R}^2$ whose support is contained in $\{\|\mathbf{x}\|\geq r\}$ for some $r>0$. This induces two bounded and continuous functions on $\mathbb{R}^2$ defined as
\[f_1(\mathbf{x})=\frac{\mathrm{dist}(\mathbf{x},U_1)}{\mathrm{dist}(\mathbf{x},U_1)+
\mathrm{dist}(\mathbf{x},U_2)}f(\mathbf{x})\;\;\;\;\;\mathrm{and}\;\;\;\;\;f_2(\mathbf{x})=
\frac{\mathrm{dist}(\mathbf{x},U_2)}{\mathrm{dist}(\mathbf{x},U_1)+
\mathrm{dist}(\mathbf{x},U_2)}f(\mathbf{x})\] for $\mathbf{x}\in(U_1\cap U_2)^c$, and $f_1(\mathbf{x})=0=f_2(\mathbf{x})$ for $\mathbf{x}\in U_1\cap U_2$, where $U_1=\{\mathbf{x}:|\mathbf{x}^{(1)}|\leq r/2\}$ and $U_2=\{\mathbf{x}:|\mathbf{x}^{(2)}|\leq r/2\}$. Clearly, $f=f_1+f_2$, and the supports of $f_1$
and $f_2$ are bounded away from the $s$- and $t$-axis, respectively. Then the weak convergence of $\{\sigma_{1n}\}$ and $\{\sigma_{2n}\}$ yields that
\begin{align*}
\lim_{n\to\infty}\int_{\mathbb{R}^2}f(s,t)\;d\tau_n(s,t)&=\lim_{n\to\infty}
\left(\int_{\mathbb{R}^2}f_1(s,t)\;d\tau_n(s,t)
+\int_{\mathbb{R}^2}f_2(s,t)\;d\tau_n(s,t)\right) \\
&=\lim_{n\to\infty}\left(\int_{\mathbb{R}^2}f_1(s,t)\frac{1+s^2}{s^2}\;d\sigma_{1n}(s,t)
+\int_{\mathbb{R}^2}f_2(s,t)\frac{1+t^2}{t^2}\;d\sigma_{2n}(s,t)\right) \\
&=\int_{\mathbb{R}^2}f_1(s,t)\frac{1+s^2}{s^2}\;d\sigma_1(s,t)
+\int_{\mathbb{R}^2}f_2(s,t)\frac{1+t^2}{t^2}\;d\sigma_2(s,t) \\
&=\int_{\mathbb{R}^2}f\;d\tau,
\end{align*} which verifies $\tau_n\Rightarrow_\mathbf{0}\tau$.

To verify the statement (IV), it suffices to prove the existence of the following limits and the equalities:
\begin{equation} \label{limitCD1}
a:=\lim_{\epsilon\to0^+}\limsup_{n\to\infty}\int_{\|\mathbf{x}\|<\epsilon}s^2\;d\tau_n(s,t)=
\lim_{\epsilon\to0^+}\liminf_{n\to\infty}\int_{\|\mathbf{x}\|<\epsilon}s^2\;d\tau_n(s,t),
\end{equation}
\begin{equation} \label{limitCD2}
b:=\lim_{\epsilon\to0^+}\limsup_{n\to\infty}\int_{\|\mathbf{x}\|<\epsilon}t^2\;d\tau_n(s,t)=
\lim_{\epsilon\to0^+}\liminf_{n\to\infty}\int_{\|\mathbf{x}\|<\epsilon}t^2\;d\tau_n(s,t),
\end{equation} and
\begin{equation} \label{limitCD3}
c:=\lim_{\epsilon\to0^+}\limsup_{n\to\infty}\int_{\|\mathbf{x}\|<\epsilon}st\;d\tau_n(s,t)=
\lim_{\epsilon\to0^+}\liminf_{n\to\infty}\int_{\|\mathbf{x}\|<\epsilon}st\;d\tau_n(s,t).
\end{equation} First of all, the limits
\[\lim_{\epsilon\to0^+}\limsup_{n\to\infty}\int_{\|\mathbf{x}\|<\epsilon}\frac{s^2}{1+s^2}\;d\tau_n(s,t)=
\lim_{\epsilon\to0^+}\liminf_{n\to\infty}\int_{\|\mathbf{x}\|<\epsilon}\frac{s^2}{1+s^2}\;d\tau_n(s,t)\]
exist and equal $\sigma_1(\{\mathbf{0}\})$ by the weak convergence of $\sigma_{1n}$ to $\sigma_1$. For any $\epsilon>0$, picking an $\epsilon'\in[\epsilon,2\epsilon]$ so that
$\sigma_1(\{\|\mathbf{x}\|=\epsilon'\})=0$ (such an $\epsilon'$ exists because $\sigma_1$ is a finite measure), we then have
\begin{align*}
\limsup_{n\to\infty}\int_{\|\mathbf{x}\|<\epsilon}\left(s^2-\frac{s^2}{1+s^2}\right)d\tau_n(s,t)
&=\limsup_{n\to\infty}\int_{\|\mathbf{x}\|<\epsilon}s^2\;d\sigma_{1n}(s,t) \\
&\leq\limsup_{n\to\infty}\int_{\|\mathbf{x}\|<\epsilon'}s^2\;d\sigma_{1n}(s,t) \\
&=\int_{\|\mathbf{x}\|<\epsilon'}s^2\;d\sigma_1(s,t) \leq\epsilon'^2\sigma_1(\mathbb{R}^2) \to0\;\;\mathrm{as}\;\;\epsilon\to0^+.
\end{align*} Hence (\ref{limitCD1}) holds and $a=\sigma_1(\{\mathbf{0}\})$. Similarly, (\ref{limitCD2}) holds true and $b=\sigma_2(\{\mathbf{0}\})$.

One can also show that the existence of the limits in (\ref{limitCD3}) is equivalent to that of the limits
\begin{equation} \label{limitCD3'}
\lim_{\epsilon\to0^+}\limsup_{n\to\infty}\int_{\|\mathbf{x}\|<\epsilon}\frac{st\;d\tau_n(s,t)}{(1+s^2)(1+t^2)}=
\lim_{\epsilon\to0^+}\liminf_{n\to\infty}\int_{\|\mathbf{x}\|<\epsilon}\frac{st\;d\tau_n(s,t)}{(1+s^2)(1+t^2)},
\end{equation} and all limits are the same if they exist. Next notice that $1\wedge\|\mathbf{x}\|^2\in L^1(\tau)$ according to the definition of $\tau$, and therefore
\begin{equation} \label{limitCD4}
\int_{\mathbb{R}^2}\frac{|st|}{(1+s^2)(1+t^2)}\;d\tau(s,t)\leq\frac{1}{2}
\int_{\mathbb{R}^2}\frac{s^2+t^2}{(1+s^2)(1+t^2)}\;d\tau(s,t)<\infty.
\end{equation}
We now show that the limits in (\ref{limitCD3'}) do exist and that the relation of $c$ and $\gamma$ in (\ref{quantity}) holds.
For any $\epsilon>0$, choose an $\epsilon'\in(\epsilon,2\epsilon]$ with the property $\tau(\{\|\mathbf{x}\|=\epsilon'\})=0$ (such an $\epsilon'$ exists because $\|\mathbf{x}\|^2\chi_{\{\|\mathbf{x}\|\leq1\}}\tau$ is a finite positive measure).
Consider the difference
\[I_n(\epsilon):=\int_{\mathbb{R}^2}\frac{st\;d\tau_n(s,t)}{(1+s^2)(1+t^2)}
-\int_{\mathbb{R}^2}\frac{st\;d\tau(s,t)}{(1+s^2)(1+t^2)}-
\int_{\{\|\mathbf{x}\|<\epsilon\}}\frac{st\;d\tau_n(s,t)}{(1+s^2)(1+t^2)}\] and decompose it into the sum of $J_{1n}(\epsilon)$ and $J_2(\epsilon)$, where
\[J_{1n}(\epsilon)=\int_{\{\|\mathbf{x}\|\geq\epsilon\}}\frac{st}{(1+s^2)(1+t^2)}\;d\tau_n(s,t)-
\int_{\{\|\mathbf{x}\|\geq\epsilon\}}\frac{st}{(1+s^2)(1+t^2)}\;d\tau(s,t)\]
and
\[J_2(\epsilon)=-\int_{\{\|\mathbf{x}\|<\epsilon\}}\frac{st}{(1+s^2)(1+t^2)}\;d\tau(s,t).\]
Denoting by $C_\epsilon$ the closed set $\{\epsilon\leq\|\mathbf{x}\|\leq\epsilon'\}$, we have
\begin{align*}
\lim_{\epsilon\to0^+}\limsup_{n\to\infty}
\left|\int_{C_\epsilon}\frac{st}{(1+s^2)(1+t^2)}\;d\tau_n\right|
\leq&\lim_{\epsilon\to0^+}\limsup_{n\to\infty}
\int_{C_\epsilon}\frac{|st|}{\sqrt{1+s^2}\sqrt{1+t^2}}\;d\tau_n \\
\leq&\lim_{\epsilon\to0^+}\limsup_{n\to\infty}
\left(\int_{C_\epsilon}\frac{s^2}{1+s^2}\;d\tau_n\right)^{1/2}
\left(\int_{C_\epsilon}\frac{t^2}{1+t^2}\;d\tau_n\right)^{1/2} \\
=&\lim_{\epsilon\to0^+}\limsup_{n\to\infty}
\sigma_{1n}(C_\epsilon)^{1/2}\sigma_{2n}(C_\epsilon)^{1/2} \\
\leq&\lim_{\epsilon\to0^+}\sigma_1(C_\epsilon)^{1/2}
\sigma_2(C_\epsilon)^{1/2}=0,
\end{align*} where the Cauchy-Schwarz inequality was used in the second inequality and the assumption that $\sigma_{jn}\Rightarrow\sigma_j$ for $j=1,2$ was used in the last one. Moreover, by Proposition \ref{Portmanteau} we have
\[\lim_{n\to\infty}\int_{\{\|\mathbf{x}\|>\epsilon'\}}\frac{st}{(1+s^2)(1+t^2)}\;d\tau_n(s,t)=
\int_{\{\|\mathbf{x}\|>\epsilon'\}}\frac{st}{(1+s^2)(1+t^2)}\;d\tau(s,t).\]
We now can conclude from the above discussions and (\ref{limitCD4}) that $\lim_{\epsilon\to0^+}\limsup_{n\to\infty}|J_{1n}(\epsilon)|=0$, whence
$\lim_{\epsilon\to0^+}\limsup_{n\to\infty}|I_n(\epsilon)|=0$ by (\ref{limitCD4}) again.
Finally, the assumption that the first integral in $I_n(\epsilon)$ converges to $\gamma$ as $n\to\infty$ yields that
the limits in (\ref{limitCD3'}) exist and equal, and the relation (\ref{quantity}) holds. Hence (\ref{limsupinf}) is proved.

The positive semi-definiteness of the matrix $\mathbf{A}$ is an easy application of the Cauchy-Schwarz inequality to (\ref{limitCD1}), (\ref{limitCD2}) and (\ref{limitCD3}). It is easy to verify that the limits in \eqref{IV} define a non-negative quadratic form $\langle\mathbf{Au},\mathbf{u}\rangle$ for any $\mathbf{u}\in\mathbb{R}^2$.

Conversely, suppose that \eqref{III} and \eqref{IV} in Condition \ref{cond2} hold. Denote by $Q(\mathbf{u})$ the finite quantity in \eqref{IV} for any $\mathbf{u}\in\mathbb{R}^2$, and define positive planar measures $\sigma_1$ and $\sigma_2$ as
\[\sigma_1=\frac{s^2}{1+s^2}\;\tau+Q((1,0))\delta_{\mathbf{0}}\;\;\;\;\;\mathrm{and}\;\;\;\;\;
\sigma_2=\frac{t^2}{1+t^2}\;\tau+Q((0,1))\delta_{\mathbf{0}}.\]
Note that measures $\sigma_1$ and $\sigma_2$ are both finite. To see this, it suffices to show that $\chi_{\{\|\mathbf{x}\|\leq1\}}\|\mathbf{x}\|^2\in L^1(\tau)$. Take a sequence $\{\epsilon_k\}_{k\geq1}$ such that $\epsilon_k\downarrow 0$ as $k\to \infty$ and $\tau(\{\|\mathbf{x}\|=\epsilon_k\})=0$ for each $k$. Then condition \eqref{IV} shows that for all $k<j$ large enough, one has
\[\limsup_{n\to\infty}\int_{\{\epsilon_j<\|\mathbf{x}\|<\epsilon_k\}}s^2\;d\tau_n(s,t)\leq Q((1,0))+1,\] which gives
\[\int_{\{\epsilon_j<\|\mathbf{x}\|<\epsilon_k\}}s^2\;d\tau(s,t)\leq Q((1,0))+1\]
by Proposition \ref{Portmanteau}. Letting $j\to\infty$ allows us to obtain that
\[\int_{\{\|\mathbf{x}\|<\epsilon_k\}}s^2\;d\tau(s,t)\leq Q((1,0))+1\] by monotone convergence theorem. This shows that $\chi_{\{\|\mathbf{x}\|\leq1\}}s^2\in L^1(\tau)$, as well as $\chi_{\{\|\mathbf{x}\|\leq1\}}t^2\in L^1(\tau)$, as desired. Now we are ready to prove that $\sigma_{1n}$ converges to $\sigma_1$ weakly. Let $f$ be a bounded and continuous function on $\mathbb{R}^2$. Then we have the following estimate
\begin{align*}
&\left|\int_{\mathbb{R}^2}f(\mathbf{x})\;d\sigma_{1n}(\mathbf{x})-
\int_{\mathbb{R}^2}f(\mathbf{x})\;d\sigma_1(\mathbf{x})\right| \\
\leq&\int_{\{\|\mathbf{x}\|<\epsilon_k\}}|f(\mathbf{x})-f(\mathbf{0})|\;d\sigma_{1n}(\mathbf{x})+
|f(\mathbf{0})||\sigma_{1n}(\{\|\mathbf{x}\|<\epsilon_k\})-Q((1,0))| \\
+&\int_{\{0<\|\mathbf{x}\|<\epsilon_k\}}|f(\mathbf{x})|\;d\sigma_1(\mathbf{x})
+\left|\int_{\{\|\mathbf{x}\|\geq\epsilon_k\}}f(\mathbf{x})\;d\sigma_{1n}(\mathbf{x})
-\int_{\{\|\mathbf{x}\|\geq\epsilon_k\}}f(\mathbf{x})\;d\sigma_1(\mathbf{x})\right| \\
=:&I_{1n}(k)+I_{2n}(k)+I_{3}(k)+I_{4n}(k),
\end{align*}
First choosing $\mathbf{u}=(1,0)$ in (\ref{limsupinf}) gives that
\[\lim_{k\to\infty}\limsup_{n\to\infty}I_{1n}(k)\leq
\lim_{k\to\infty}\left[\max_{\|\mathbf{x}\|<\epsilon_k}|f(\mathbf{x})-f(\mathbf{0})|\right]
\cdot\left(\limsup_{n\to\infty}
\int_{\{\|\mathbf{x}\|<\epsilon_k\}}s^2\;d\tau_n(s,t)\right)=0.\]
Since
\begin{align*}
\lim_{k\to\infty}\limsup_{n\to\infty}\left|\int_{\{\|\mathbf{x}\|<\epsilon_k\}}\frac{s^2}{1+s^2}\;d\tau_n-
\int_{\{\|\mathbf{x}\|<\epsilon_k\}}s^2\;d\tau_n\right|
\leq\lim_{k\to\infty}\epsilon_k^2\cdot\left(\limsup_{n\to\infty}\int_{\{\|\mathbf{x}\|<\epsilon_k\}}
s^2\;d\tau_n\right)=0,
\end{align*} it follows that
\[\lim_{k\to\infty}\limsup_{n\to\infty}\big|\sigma_{1n}(\{\|\mathbf{x}\|<\epsilon_k\})-Q((1,0))\big|=0,\] whence
$\lim_{k\to\infty}\limsup_{n\to\infty}I_{2n}(k)=0$.
We also have
$\limsup_{k\to\infty}I_{3}(k)\leq\|f\|_\infty\lim_{k\to\infty}\sigma_1(\{0<\|\mathbf{x}\|<\epsilon_k\})=0$ and $\limsup_{n\to\infty}I_{4n}(k)=0$ because $\sigma_1$ is a finite measure and $\tau_n\Rightarrow_\mathbf{0}\tau$.
Thus we have shown the weak convergence of the sequence $\{\sigma_{1n}\}$ in condition \eqref{I}. Similarly, $\sigma_{2n}\Rightarrow\sigma_2$.

Finally, we decompose the desired integral in condition \eqref{II} into the sum
\[\int_{\|\mathbf{x}\|<\epsilon_k}\frac{st}{(1+s^2)(1+t^2)}\;d\tau_n(s,t)+
\int_{\|\mathbf{x}\|\geq\epsilon_k}\frac{st}{(1+s^2)(1+t^2)}\;d\tau_n(s,t).\]
As $n\to\infty$ and $k\to\infty$, the first integral tends to $[Q((1,1))-Q((1,0))-Q((0,1))]/2$, while the second integral tends to
\[\int_{\mathbb{R}^2}\frac{st}{(1+s^2)(1+t^2)}\;d\tau(s,t)\] by Proposition \ref{Portmanteau} and the fact that $1\wedge\|\mathbf{x}\|^2 \in L^1(\tau)$. Hence condition \eqref{II} is verified and the proof is complete.
\end{pf} \qed

We are in a position to prove the equivalence between classical and bi-free limit theorems for non-identical distributions.

\begin{thm} \label{limitthm1}
Let $\{k_n\}_{n=1}^\infty\subset\mathbb{N}$ be strictly increasing, $\{\mu_{nk}\}_{1\leq n,1\leq k\leq k_n}\subset\mathscr{P}_{\mathbb{R}^2}$ be an infinitesimal triangular array and $\{\mathbf{v}_n\}_{n=1}^\infty\subset\mathbb{R}^2$. With the notations in \eqref{vnk}, \eqref{barmunk}, and \eqref{nun}, the following are equivalent.
\begin{enumerate} [$\qquad(1)$]
\item\label{Assertion1}
The sequence
\[\mu_{n1}*\mu_{n2}*\cdots*\mu_{nk_n}*\delta_{\mathbf{v}_n}\]
converges weakly to some probability measure $\mu_*$ on $\mathbb{R}^2$.

\item\label{Assertion2} The sequence
\[\mu_{n1}\bfconv\mu_{n2}\bfconv\cdots\bfconv\mu_{nk_n}\bfconv\delta_{\mathbf{v}_n}\]
converges weakly to some probability measure $\mu_{\bfconvv}$ on $\mathbb{R}^2$.

\item\label{Assertion3} {\normalfont Condition \ref{cond1}}\eqref{I} and \eqref{II} hold, and the vector in
\eqref{vector} converges to some vector $\mathbf{v} \in \mathbb{R}^2$.

\item\label{Assertion4} {\normalfont Condition \ref{cond2}}\eqref{III} and \eqref{IV} hold, and the vector in
\eqref{vector} converges to some vector $\mathbf{v} \in \mathbb{R}^2$.
\end{enumerate}
If assertions \eqref{Assertion1}-\eqref{Assertion4} hold, then $\mu_*$ and $\mu_{\bfconvv}$ are $*$-infinitely divisible and $\bfconvv$-infinitely divisible distributions with $(\mathbf{v},\mathbf{A},\tau)$ as the classical and bi-free characteristic triplet, respectively, where $\mathbf{A}$ is defined as in \eqref{matrix}.
\end{thm}

\begin{pf} The equivalences $\eqref{Assertion2}\Leftrightarrow\eqref{Assertion3}$ and $\eqref{Assertion3}\Leftrightarrow\eqref{Assertion4}$ were already respectively proved in Theorem \ref{nonidentical} and Lemma \ref{mainlem}, while the equivalence $\eqref{Assertion1}\Leftrightarrow\eqref{Assertion4}$ can be obtained by \cite[Theorem 3.2.2 and (3.52),(3.53),(3.54)]{limitthm}. We remark here that some results cited from \cite{limitthm} contain errors, and the reader is referred to the list of errata of the book put on the webpage of one of the authors.

It remains to show that $\phi_{\mu_{\bfconvv}}$ extends analytically to $(\mathbb{C}\backslash\mathbb{R})^2$ and admits an integral representation of the form (\ref{BFLK}). For any $\epsilon>0$ and $(z,w)\in(\mathbb{C}\backslash\mathbb{R})^2$, let
\[\mathcal{P}_n(z,w,\epsilon)=\int_{\{\|\mathbf{x}\|\geq\epsilon\}}
\left[\frac{zw}{(z-s)(w-t)}-1-\frac{z^{-1}s+w^{-1}t}{1+s^2+t^2}\right]d\tau_n(s,t)\] and
\[\mathcal{G}_n(z,w,\epsilon)=\int_{\{\|\mathbf{x}\|<\epsilon\}}
\left[\frac{zw}{(z-s)(w-t)}-1-\frac{z^{-1}s+w^{-1}t}{1+s^2+t^2}\right]d\tau_n(s,t).\] Notice that the integrand in the integral can be rewritten as
\[\frac{1}{z(z-s)}\frac{s^2}{1+s^2+t^2}+\frac{1}{w(w-t)}\frac{t^2}{1+s^2+t^2}+\frac{st}{(z-s)(w-t)}.\] Then
choosing $\epsilon$ so that $\tau(\{\|\mathbf{x}\|=\epsilon\})=0$ shows that
\[\mathcal{P}(z,w):=\lim_{\epsilon\to0^+}\lim_{n\to\infty}\mathcal{P}_n(z,w,\epsilon)=\int_{\mathbb{R}^2}
\left[\frac{zw}{(z-s)(w-t)}-1-\frac{z^{-1}s+w^{-1}t}{1+s^2+t^2}\right]
d\tau(s,t),\] where we used Proposition \ref{Portmanteau} and the fact that $1\wedge\|\mathbf{x}\|^2\in L^1(\tau)$. On the other hand, (\ref{limitCD1}), (\ref{limitCD2}) and (\ref{limitCD3}) yield that
\[\mathcal{G}(z,w):=\lim_{\epsilon\to0^+}\lim_{n\to\infty}\mathcal{G}_n(z,w,\epsilon)
=\frac{\sigma_1(\{\mathbf{0}\})}{z^2}+
\frac{c}{zw}+\frac{\sigma_2(\{\mathbf{0}\})}{w^2}.\] Then according to (\ref{BFLKyv}), $\phi_{\mu_{\bfconvv}}(z,w)$ agrees with
\[\frac{\mathbf{v}^{(1)}}{z}+\frac{\mathbf{v}^{(2)}}{w}+\mathcal{G}(z,w)+\mathcal{P}(z,w)\] when $z=iy$ and $w=iv$ with $|y|$ and $|v|$ large, and hence they agree on $(\mathbb{C}\backslash\mathbb{R})^2$ by analytic extension. The last assertion regarding $\mu_{\bfconvv}$ follows from Theorem \ref{bifreeID}. This finishes the proof.
\end{pf} \qed

The classical and bi-free L\'evy-Hin\v{c}in representations \eqref{LK1} and \eqref{BFLK} establish a bijective relation $\Lambda$ between the sets $\mathcal{ID}(*)$ and $\mathcal{ID}(\bfconvv)$:
\begin{equation} \label{bijection}
\Lambda\big(\mu_*^{(\mathbf{v,A},\tau)}\big)=\mu_{\bfconvv}^{(\mathbf{v,A},\tau)}
\end{equation}
for any infinitely divisible law $\mu_*^{(\mathbf{v,A},\tau)}$ with classical characteristic triplet $(\mathbf{v,A},\tau)$. Under this bijection, classical Gaussian and (compound) Poisson
distributions are respectively mapped to bi-free Gaussian and bi-free (compound) Poisson
distributions (see Example \ref{example}). Furthermore, Theorem \ref{limitthm1} and this bijection establish a transfer principle for limit theorems.

The limit theorem for the identically distributed random variables is formulated below.

\begin{thm} \label{limitthm2}
Let $\{\mu_n\}_{n=1}^\infty$ be a sequence in $\mathscr{P}_{\mathbb{R}^2}$ and let $\{k_n\}_{n=1}^\infty\subset\mathbb{N}$ be strictly increasing. Then the statements \eqref{Cconv}--\eqref{III-IV} are equivalent.
\begin{enumerate} [$\qquad(1)$]
\item\label{Cconv} The measure $\mu_n^{*k_n}$ converges weakly
to some probability measure $\mu_*$ on $\mathbb{R}^2$.

\item\label{BFconv} The sequence $\mu_n^{\bfconvv k_n}$ converges
weakly to some probability measure $\mu_{\bfconvv}$ on
$\mathbb{R}^2$.

\item\label{I-II} {\normalfont Condition \ref{cond1}}\eqref{I} and \eqref{II} hold with $\tau_n=k_n\mu_n$,
and the limit
\begin{equation} \label{vector2}
\lim_{n\to\infty}k_n\int_{\mathbb{R}^2}\frac{\mathbf{x}}{1+\|\mathbf{x}\|^2}\;d\mu_n(\mathbf{x})=\mathbf{v}
\end{equation}
exists.

\item\label{III-IV} {\normalfont Condition \ref{cond2}}\eqref{III} and \eqref{IV} hold with $\tau_n=k_n\mu_n$,
and the limit in \eqref{vector2} exists.
\end{enumerate}
If assertions \eqref{Cconv} through \eqref{III-IV} hold, then $\mu_*$ and $\mu_{\bfconvv}$ are $*$-infinitely divisible and $\bfconvv$-infinitely divisible distributions with $(\mathbf{v},\mathbf{A},\tau)$ as the classical and bi-free characteristic triplet, respectively, where $\mathbf{A}$ is defined as in \eqref{matrix}.
\end{thm}

\begin{remark}{\normalfont
Due to the recent work of Gu and Skoufranis \cite[Theorem 5.12]{GS}, the above conditions \eqref{Cconv}--\eqref{III-IV} are further equivalent to the statement that the sequence $\mu_n^{\uplus\uplus k_n}$ converges
weakly to some probability measure $\mu_{\uplus\uplus}$ on $\mathbb{R}^2$, where $\uplus\uplus$ is the bi-boolean convolution. The limit distribution $\mu_{\uplus\uplus}$ may also be characterized by a bi-boolean characteristic triplet.}
\end{remark}

\begin{pf} The equivalence $\eqref{Cconv}\Leftrightarrow\eqref{BFconv}$ will follow immediately by choosing $\mu_{nk}=\mu_n$ for $k=1,\dots,k_n$ and $\mathbf{v}_n=\mathbf{0}$ in Theorem \ref{limitthm1} once the infinitesimality of $\{\mu_n\}$ is verified in (1) and (2). In (1), we have $\mu_n^{(j)}\Rightarrow\delta_0$ for $j=1,2$ by \cite{Kolmogorov} (\S14, Theorem 4), whence $\mu_n\Rightarrow\delta_\mathbf{0}$ by \eqref{eqA}. In (2), we see that $\phi_{\mu_n}\to0=\phi_{\delta_\mathbf{0}}$ uniformly on compact sets of $\Gamma^2$ and $\phi_{\mu_n}(z,w)=o(1)$ uniformly in $n$ as $z,w\to\infty$ with $(z,w)\in\Gamma^2$ by Proposition \ref{weakconvphi}, whence $\mu_n\Rightarrow\delta_\mathbf{0}$ by Proposition \ref{weakconvphi} again. Hence the infinitesimality is verified in both situations.

The equivalence $\eqref{BFconv}\Leftrightarrow\eqref{I-II}$ was already proved in \cite[Theorem 3.2]{HW}, while the equivalence $\eqref{I-II}\Leftrightarrow\eqref{III-IV}$ can be obtained by applying Lemma \ref{mainlem} to the positive measures $\tau_n=k_n\mu_n$. That the limiting distribution $\mu_{\bfconvv}$ is $\bfconvv$-infinitely divisible distributions with the desired bi-free characteristic triplet follows from the last part of the proof of Theorem \ref{limitthm1}.
\end{pf} \qed

\begin{remark} {\normalfont In the proof of \eqref{BFconv} $\Leftrightarrow$ \eqref{I-II} in Theorem \ref{limitthm2}, a bi-free limit theorem on identical distributions was employed \cite{HW}. As one might expect, a more direct proof based on Theorem \ref{limitthm1} without referring to any other type of limit theorems exists, but it is not a short one. More precisely, what one really needs to show is the equivalence of the following two statements:
\begin{enumerate}[\qquad(i)]
\item \normalfont Condition \ref{cond2} holds with $\tau_n=k_n\mu_n$
and there exists some $\mathbf{v}\in \mathbb{R}^2$ so that (\ref{vector2}) holds;
\item \normalfont Condition \ref{cond2} holds with $\tau_n = k_n \mathring{\mu}_n$ and
there exists some $\mathbf{v}\in \mathbb{R}^2$ so that
\[\lim_{n\to\infty}k_n\left[\mathbf{u}_n+\int_{\mathbb{R}^2}\frac{\mathbf{x}}{1+\|\mathbf{x}\|^2} \;d\mathring{\mu}_{n}(\mathbf{x})\right]=\mathbf{v},\]
\end{enumerate}
where $\mathring{\mu}_n$ is the shift of $\mu_n$ by the vector $\mathbf{u}_n:=\int_{\|\mathbf{x}\|<L}\mathbf{x}\,d\mu_n(\mathbf{x})$.
Some elaboration and techniques are needed to show this equivalence. We leave the proof to the interested reader.}
\end{remark}

\section{Stable laws in bi-free probability}
In this section we define and study bi-free stable distributions, and show that they arise naturally in limit theorems. The presented result establishes the coincidence of the domains of attraction in classical probability and bi-free probability.

For any $\lambda>0$, denote by $D_\lambda$ the dilation operator on measures $\rho$ on $\mathbb{R}^d$, i.e. for any Borel set $B\subset\mathbb{R}^d$,
\[(D_\lambda\rho)(B)=\rho(\{\lambda^{-1}\mathbf{x}:\mathbf{x}\in B\}).\]

\begin{pdef} {\normalfont Let $\star$ be a binary operation on the set $\mathscr{P}_{\mathbb{R}^2}$. A planar probability distribution $\mu$ is said to be $\star$-\emph{stable} if for any $a,b>0$, there exist some $c>0$ and some vector $\mathbf{u}\in\mathbb{R}^2$ so that
\[(D_a\mu)\star(D_b\mu)=(D_c\mu)\star\delta_\mathbf{u}.\]}
\end{pdef}

The classification of $*$-stable distributions is known, but the authors could not find a reference including a complete proof, so the statement with a proof is provided below. We say that a probability measure is \emph{non-trivial} if it is not a delta measure.

\begin{thm} \label{Cstable}
A non-trivial planar probability measure is $*$-stable if and only if either
\begin{enumerate} [$\qquad(1)$]
\item\label{CstableG} it is a Gaussian distribution or
\item\label{CstableA} {it is $*$-infinitely divisible and admits the $*$-characteristic triplet
$(\mathbf{v},\mathbf{0},\tau)$ with $\tau$ of the form
\[d\tau(\mathbf{x})=\frac{1}{r^{1+\alpha}}\;drd\Theta(\omega),\] where $\alpha\in(0,2)$, $\Theta$ is a finite positive measure on the unit circle
$\mathbb{T}$ and $\mathbf{x}=r\omega$ with $r>0$ and $\omega\in\mathbb{T}$.}
\end{enumerate}
\end{thm}

\begin{pf} Suppose that $\mu$ is non-trivial and $*$-stable. We may assume that the marginal law $\mu^{(1)}$ is non-trivial. Then for any $a,b>0$, there exist some $c>0$ and $\mathbf{u}\in\mathbb{R}^2$ so that
\begin{equation}\label{StableC2}
(D_a\mu^{(1)})\ast(D_b\mu^{(1)})=\big[(D_a\mu)\ast(D_b\mu)\big]^{(1)}=\big[(D_c\mu)\ast\delta_\mathbf{u}\big]^{(1)}=
(D_c\mu^{(1)})\ast\delta_{\mathbf{u}^{(1)}}.
\end{equation}
This shows that $\mu^{(1)}$ is $\ast$-stable, and hence its L\'evy measure is either a zero measure or of the form $d\rho(x):=c_1x^{-\alpha-1}\chi_{(0,\infty)}(x)dx+c_2|x|^{-\alpha-1}\chi_{(-\infty,0)}(x)dx$ for some $\alpha\in(0,2)$ and $c_1,c_2\geq0$ with $c_1+c_2>0$ (see \S34 in \cite{Kolmogorov}). If $\rho\neq0$, one can check from characteristic functions that the constant $c$ is uniquely determined by the relation $c^\alpha=a^\alpha+b^\alpha$.
In the first case, $\mu^{(1)}$ is Gaussian, and hence $c^2=a^2+b^2$ from characteristic functions again, which is realized as $\alpha=2$.

To obtain the desired result, let $\alpha\in(0,2]$ be fixed and consider the family $(\mu_\lambda)_{\lambda>0}$, where $\mu_\lambda=D_{\lambda^{1/\alpha}}\mu$. Then using the $*$-stability of $\mu$ and the relation $a^\alpha+b^\alpha=c^\alpha$ shows that
\begin{equation} \label{semigroup}
\mu_{\lambda_1+\lambda_2}=\mu_{\lambda_1}\ast\mu_{\lambda_2}\ast
\delta_{\mathbf{u}(\lambda_1,\lambda_2)}
\end{equation}
for some $\mathbf{u}(\lambda_1,\lambda_2)\in\mathbb{R}^2$, which clearly gives the $*$-infinite divisibility of $\mu$ and each $\mu_\lambda$. Let $(\mathbf{v},\mathbf{A},\tau)$ be the $*$-characteristic triplet of $\mu$. Then $\mu_\lambda$ admits the $*$-characteristic triplet $(\mathbf{v}(\lambda),\lambda^{2/\alpha}\mathbf{A},D_{\lambda^{1/\alpha}}\tau)$ for some $\mathbf{v}(\lambda)\in\mathbb{R}^2$. Moreover, (\ref{semigroup}) yields the following two relations for any $\lambda_1,\lambda_2>0$:
\begin{equation} \label{dilation1C}
(\lambda_1+\lambda_2)^{2/\alpha}\mathbf{A}=\lambda_1^{2/\alpha}\mathbf{A}+\lambda_2^{2/\alpha}\mathbf{A}
\end{equation}
and
\begin{equation} \label{dilation2C}
D_{(\lambda_1+\lambda_2)^{1/\alpha}}\tau=D_{\lambda_1^{1/\alpha}}\tau+D_{\lambda_2^{1/\alpha}}\tau.
\end{equation}

To continue the proof, let $\Omega\subset\mathbb{T}$ be a fixed Borel set. By restricting the measures appearing in (\ref{dilation2C}) on the set $\{r\Omega:r\geq1\}$, one can infer that the function $f(\lambda)=\tau(\{r\Omega:r\geq\lambda^{-1/\alpha}\})$ satisfies Cauchy's functional equation $f(\lambda_1+\lambda_2)=f(\lambda_1)+f(\lambda_2)$. Since $f$ is increasing on $(0,\infty)$, it is measurable there, and hence $f(\lambda)=\lambda f(1)$ for any $\lambda>0$. This allows us to obtain that
\[\tau(\{r\Omega:r\geq\lambda\})=\lambda^{-\alpha}\tau(\{r\Omega:r\geq1\})\] for any $\lambda>0$ and any Borel set $\Omega\subset\mathbb{T}$. Hence the finite positive measure
\[\Theta(\Omega)=\alpha\int_{[1,\infty)\times\Omega}\;d\tau(r,\omega),\;\;\;\;\;\Omega\subset\mathbb{T},\] gives us the desired one. If $\alpha=2$, then only $\tau=0$ is allowed in order to fit the condition $1\wedge\|\mathbf{x}\|^2\in L^1(\tau)$, in which case $\mu$ is a Gaussian. If $\alpha<2$, then (\ref{dilation1C}) holds for any $\lambda_1,\lambda_2>0$ if and only if $\mathbf{A}=0$. For the converse, it is clear that $\mu$ is $\ast$-stable either in the case \eqref{CstableG} or \eqref{CstableA}.
\end{pf} \qed

The $\bfconvv$-stable distributions are classified as follows.

\begin{thm} \label{Bstable} A non-trivial planar probability measure is $\bfconvv$-stable if and only if either
\begin{enumerate} [$\qquad(1)$]
\item {it is a bi-free Gaussian distribution or}
\item {it is $\bfconvv$-infinitely divisible and it has a $\bfconvv$-characteristic triplet $(\mathbf{v},\mathbf{0},\tau)$ with $\tau$ of the form
\[d\tau(\mathbf{x})=\frac{1}{r^{1+\alpha}}\;drd\Theta(\omega),\] where $\alpha\in(0,2)$, $\Theta$ is a finite positive measure on the unit circle
$\mathbb{T}$ and $\mathbf{x}=r\omega$ with $r>0$ and $\omega\in\mathbb{T}$.}
\end{enumerate}
\end{thm}

\begin{pf} Suppose that $\mu$ is non-trivial and $\bfconvv$-stable. Further suppose that $\mu^{(1)}$ is non-trivial. Then it follows from Proposition \ref{bifreeproj2} that
\[(D_a\mu^{(1)})\boxplus(D_b\mu^{(1)})=\big[(D_a\mu)\bfconv(D_b\mu)\big]^{(1)}=\big[(D_c\mu)\bfconv\delta_\mathbf{u}\big]^{(1)}=
(D_c\mu^{(1)})\boxplus\delta_{\mathbf{u}^{(1)}}.\] This gives the $\boxplus$-stability of $\mu^{(1)}$, and hence $\phi_{\mu^{(1)}}'(z)=\beta z^{-\alpha}$ for some $\alpha\in(0,2]$ and $\beta\in\mathbb{C}\backslash\{0\}$ by Lemma 7.4 and Theorem 7.5 of \cite{BerVoicu93}. Since $\phi_{D_\lambda\mu^{(1)}}(z)=\lambda\phi_{\mu^{(1)}}(z/\lambda)$ for any $\lambda>0$, one can conclude that $a$, $b$ and $c$ satisfy the relation $c^\alpha=a^\alpha+b^\alpha$.

As in the proof of Theorem \ref{Cstable}, we consider the measures $\mu_\lambda=D_{\lambda^{1/\alpha}}\mu$, $\lambda>0$. Then the $\bfconvv$-stability of $\mu$ shows that
$\mu_{\lambda_1+\lambda_2}=\mu_{\lambda_1}\bfconv\mu_{\lambda_2}\bfconv\delta_{\mathbf{u}(\lambda_1,\lambda_2)}$
for some vector $\mathbf{u}(\lambda_1,\lambda_2)\in\mathbb{R}^2$, which gives the infinite divisibility of $\mu$ and each $\mu_\lambda$. If $(\mathbf{v},\mathbf{A},\tau)$ is the bi-free characteristic triplet of $\mu$, then the identity $\phi_{D_\lambda\mu}(z,w)=\phi_\mu(z/\lambda,w/\lambda)$, which holds for $(z,w)$ in the common domain of these transforms, yields that $\mu_\lambda$ admits the bi-free characteristic triplet $(\mathbf{v}(\lambda),\lambda^{2/\alpha}\mathbf{A},D_{\lambda^{1/\alpha}}\tau)$ for some $\mathbf{v}(\lambda)\in\mathbb{R}^2$.
Then the remaining proof is similar to that of Theorem \ref{Cstable}.
\end{pf} \qed

All $\bfconvv$-stable distributions are $\bfconvv$-infinitely divisible. The number $\alpha\in(0,2]$ is called \emph{stability index} of $\mu_*$ and $\mu_{\bfconvv}$. It is shown in the proof that $\mu_*$ and $\Lambda(\mu_*)$ have the same stability index, particularly, the stability index of Gaussian and bi-free Gaussian are both two.

In \cite{Rva}, Rva\v{c}eva investigated the limiting distribution of random vectors
\begin{equation} \label{ST}
\frac{X_1+\cdots+X_n}{b_n}+\mathbf{u}_n,
\end{equation} where $\{X_n\}_{n\geq1}$ are i.i.d.\ random vectors, $b_n>0$ and $\mathbf{u}_n\in\mathbb{R}^2$. It turns out that the set of all possible limiting distributions in (\ref{ST}) equals the set of $*$-stable distributions from the arguments in \cite[\S33]{Kolmogorov}. The limit theorem of this type in the bi-free setting is considered as follows.

\begin{thm} \label{domainatt}
Let $\nu$ be a planar probability distribution, $b_n>0$ and $\mathbf{u}_n\in\mathbb{R}^2$ for $n=1,2,\ldots$. Then the following statements are equivalent.
\begin{enumerate} [$\qquad(1)$]
\item\label{St1} The measures $(D_{1/b_n}\nu^{*n})*\delta_{\mathbf{u}_n}$
converge weakly to a probability distribution $\mu_*$ on $\mathbb{R}^2$.

\item\label{St2} The measures $(D_{1/b_n}\nu^{\bfconvv n})*\delta_{\mathbf{u}_n}$
converge weakly to a probability distribution
$\mu_{\bfconvv}$ on $\mathbb{R}^2$.
\end{enumerate}
If \eqref{St1} and \eqref{St2} hold, then $\mu_*$ and $\mu_{\bfconvv}$ are $*$-stable and $\bfconvv$-stable, respectively, whose respective $*$-characteristic triplet and $\bfconvv$-characteristic triplet coincide.
\end{thm}

\begin{pf} Applying Theorem \ref{limitthm2} to the positive integers $k_n=n$ and the measures $\mu_{n}=(D_{1/b_n}\nu)*\delta_{\mathbf{u}_n/n}$ and $\mu_n=(D_{1/b_n}\nu)\bfconv\delta_{\mathbf{u}_n/n}$ in (1) and (2), respectively, yields the desired equivalence.  The last statement is a direct consequence of the mentioned results around \eqref{ST} and the fact $\mu_{\bfconvv} = \Lambda(\mu_\ast)$ established in Theorem \ref{limitthm2}.
\end{pf} \qed

\begin{pdef} {\normalfont A measure $\nu\in\mathscr{P}_{\mathbb{R}^2}$ is said to belong to the \emph{$\star$-domain of attraction} of a $\star$-stable law $\mu_{\star}$ if there exist a sequence $\{b_n\}_{n=1}^\infty$ of positive numbers and a sequence $\{\mathbf{u}_n\}_{n=1}^\infty$ of vectors in $\mathbb{R}^2$ so that $(D_{1/b_n}\nu^{\star n})\star\delta_{\mathbf{u}_n}\Rightarrow\mu_{\star}$. Denote by $\mathbf{D}_\star(\mu_{\star})$ the $\star$-domain of attraction of a give $\star$-stable law $\mu_{\star}$.}
\end{pdef}

The $*$-domain of attraction was studied in great detail in \cite{Rva}. One can immediately conclude the following result from Theorem \ref{domainatt}.

\begin{cor} For any $*$-stable law $\mu_{\ast}$ on $\mathbb{R}^2$, $\mathbf{D}_*(\mu_\ast)=\mathbf{D}_{\bfconvv}(\Lambda(\mu_\ast))$.
\end{cor}

\section{Full distributions}
We will discuss in this section the concept of fullness which regards the supports of probability distributions introduced below:

\begin{pdef}
\emph{A Borel measure $\rho$ on $\mathbb{R}^2$ is said to be \emph{full} if it is not supported on a straight line, while $\rho$ is called $\mathcal{M}_{\mathbb{R}^2}^\mathbf{0}$-\emph{full} if it is in $\mathcal{M}_{\mathbb{R}^2}^\mathbf{0}$ and not supported on a line through the origin.}
\end{pdef}

A bivariate normal distribution is full if and only if its symmetric covariance matrix is strictly positive definite, in which case the distribution has a density. If the covariance matrix is not of full rank, then the bivariate normal distribution is non-full and does not have a density.

In the following we relate the fullness of measures in $\mathscr{P}_{\mathbb{R}^2}$ to their Cauchy transforms and bi-free $\phi$-transforms.

\begin{lem} \label{fulllem}
A measure $\mu\in\mathscr{P}_{\mathbb{R}^2}$ is non-full if and only if there exist $\alpha,\beta,\gamma\in\mathbb{R}$ so that
\begin{equation} \label{fullcond1}
(\alpha z+\beta w+\gamma)G_\mu(z,w)=\beta G_{\mu^{(1)}}(z)+\alpha G_{\mu^{(2)}}(w)
\end{equation} holds for any $(z,w)\in(\mathbb{C}\backslash\mathbb{R})^2$, in which case it is supported on the line $\alpha s+\beta t+\gamma=0$.
\end{lem}

\begin{pf} First notice that for $\alpha,\beta,\gamma\in\mathbb{R}$ and $(z,w)\in(\mathbb{C}\backslash\mathbb{R})^2$, we have
\begin{align*}
G(z,w):&=\int_{\mathbb{R}^2}\frac{\alpha s+\beta t+\gamma}{(z-s)(w-t)}\;d\mu(s,t) \\
&=(\alpha z+\beta w+\gamma)G_\mu(z,w)-\beta G_{\mu^{(1)}}(z)-\alpha G_{\mu^{(2)}}(w).
\end{align*} This clearly gives (\ref{fullcond1}) if $\mu$ is supported on $\alpha s+\beta t+\gamma=0$. Conversely, suppose that $G(z,w)=0$ holds true for $(z,w)\in(\mathbb{C}\backslash\mathbb{R})^2$. Then a simple computation shows that
\begin{equation} \label{vanish1}
\int_{\mathbb{R}^2}\frac{\alpha s+\beta t+\gamma}{(s^2+1)(t^2+1)}\;d\mu(s,t)=-\frac{\Re\big[G(i,i)-G(-i,i)\big]}{2}=0.
\end{equation}
On the other hand, considering the function $H(z,w)=zG(z,w)$ yields that
\begin{equation} \label{vanish2}
\int_{\mathbb{R}^2}\frac{s(\alpha s+\beta t+\gamma)}{(s^2+1)(t^2+1)}\;d\mu(s,t)=-\frac{\Re\big[H(i,i)-H(-i,i)\big]}{2}=0.
\end{equation}
Similarly, one can obtain that
\begin{equation} \label{vanish3}
\int_{\mathbb{R}^2}\frac{t(\alpha s+\beta t+\gamma)}{(s^2+1)(t^2+1)}\;d\mu(s,t)=0.
\end{equation} Multiplying (\ref{vanish1}), (\ref{vanish2}) and (\ref{vanish3}) by $\gamma$, $\alpha$ and $\beta$, respectively, and then adding them all together shows that
\[\int_{\mathbb{R}^2}\frac{(\alpha s+\beta t+\gamma)^2}{(s^2+1)(t^2+1)}\;d\mu(s,t)=0.\] Since $\mu$ is positive, this clearly shows that it is supported on $\alpha s+\beta t+\gamma=0$, as desired.
\end{pf} \qed

\begin{prop} \label{fullcond}
A measure $\mu\in\mathscr{P}_{\mathbb{R}^2}$ is non-full if and only if there exist $\alpha,\beta,\gamma\in\mathbb{R}$ so that
\begin{equation} \label{fullcond2}
zw(\alpha z+\beta w)\phi_\mu(z,w)=\beta w^2\phi_{\mu^{(1)}}(z)+\alpha z^2\phi_{\mu^{(2)}}(w)-\gamma zw
\end{equation}
holds for $(z,w)\in\Gamma^2$, in which case $\mu$ is supported on the line $\alpha s+\beta t+\gamma=0$.
\end{prop}

\begin{pf} With the help of Lemma \ref{fulllem}, we see that $\mu$ is supported on the line $\alpha s+\beta t+\gamma=0$ if and only if (\ref{fullcond1}) holds for $(z,w)\in(\mathbb{C}\backslash\mathbb{R})^2$ or, equivalently, the identity
\begin{equation} \label{fullcond3}
\big(\alpha F_{\mu^{(1)}}^{-1}(z)+\beta F_{\mu^{(2)}}^{-1}(w)+\gamma\big)G_\mu\left(F_{\mu^{(1)}}^{-1}(z),F_{\mu^{(2)}}^{-1}(w)\right)=\frac{\beta}{z}+\frac{\alpha}{w}
\end{equation} is valid for $(z,w)\in\Gamma^2$. Note that the function $G_\mu(F_{\mu^{(1)}}^{-1},F_{\mu^{(2)}}^{-1})$ never vanishes on $\Gamma^2$ by shrinking the domain if necessary. Then we see that $\mu$ is supported on $\alpha s+\beta t+\gamma=0$ if and only if
\begin{equation} \label{fullcond4}
(\alpha z+\beta w)\left[1-\frac{1}{zwG_\mu\left(F_{\mu^{(1)}}^{-1}(z),F_{\mu^{(2)}}^{-1}(w)\right)}\right]=-\left[\alpha\phi_{\mu^{(1)}}(z)
+\beta\phi_{\mu^{(2)}}(w)+\gamma\right]
\end{equation}
holds true for $(z,w)\in\Gamma^2$. Apparently, (\ref{fullcond2}) and (\ref{fullcond4}) are equivalent, concluding the proof.
\end{pf} \qed

\begin{thm} \label{fullMfull2}
Let $\mu$ be a $\boxplus\boxplus$-infinitely divisible distribution on $\mathbb{R}^2$ with bi-free characteristic triplet $[\mathbf{v},\mathbf{A},\tau]$. Then $\mu$ is non-full if and only if $\mathbf{A}$ is singular and $\tau$ is supported on
$\langle\mathbf{u},(s,t)\rangle=0$ with some $\mathbf{u}\neq\mathbf{0}$ in the kernel of $\mathbf{A}$, in which case $\mu$ is supported on $\langle\mathbf{u},(s,t)\rangle=\langle\mathbf{u},\mathbf{v}\rangle$.
\end{thm}

\begin{pf} We shall use Proposition \ref{fullcond} to conclude the proof. First notice that for any real numbers $\alpha$, $\beta$ and $\gamma$, the function
\begin{equation} \label{fullcond7}
\frac{\alpha z+\beta w}{zw}\phi_\mu(z,w)-\frac{\beta}{z^2}\phi_{\mu^{(1)}}(z)-\frac{\alpha}{w^2}\phi_{\mu^{(2)}}(w)+\frac{\gamma}{zw}
\end{equation}
can be expressed as
\begin{equation} \label{fullcond8}
\frac{\gamma'}{zw}+\frac{\alpha a+\beta c}{z^2w}+\frac{\alpha c+\beta b}{zw^2}
+\int_{\mathbb{R}^2}\left[\frac{\alpha s+\beta t}{(z-s)(w-t)}-
\frac{1}{zw}\frac{\alpha s+\beta t}{1+s^2+t^2}\right]d\tau(s,t)
\end{equation}
by (\ref{BFLK}), (\ref{BFLK1}) and (\ref{BFLK2}),
where $\mathbf{A}$ is as in (\ref{A1}) and
$\gamma'=\alpha\mathbf{v}^{(1)}+\beta\mathbf{v}^{(2)}+\gamma$.

Let $\mathbf{u}=(\alpha,\beta)$ and $\gamma = -\langle \mathbf{u}, \mathbf{v}\rangle$. In this case $\gamma'=0$. If $\mu$ is supported on $\langle\mathbf{u},(s,t)\rangle=\langle \mathbf{u}, \mathbf{v}\rangle$, then Proposition \ref{fullcond} yields that the function in (\ref{fullcond8}) vanishes on $(\mathbb{C}\backslash\mathbb{R})^2$. Using the technique employed in Lemma \ref{fulllem} we can obtain that
\[\int_{\mathbb{R}^2}\frac{(\alpha s+\beta t)^2}{(1+s^2)(1+t^2)}\;d\tau(s,t)=-(\alpha^2a+2\alpha\beta c+\beta^2b)
=-\langle\mathbf{Au},\mathbf{u}\rangle.\] Since $\mathbf{A}\geq0$, it follows that $\mathbf{Au}=0$ and $\tau$ is supported on the line $\alpha s+\beta t=0$. Conversely, if $\mathbf{Au}=0$ and $\tau$ is supported on the line $\alpha s+\beta t=0$, then using Proposition \ref{fullcond} and (\ref{fullcond8}) again shows that $\mu$ is supported on $\alpha s+\beta t=\alpha\mathbf{v}^{(1)}+\beta\mathbf{v}^{(2)}$, as desired.
\end{pf} \qed

The following results are both direct consequences of Theorem \ref{fullMfull2}.

\begin{cor} A bi-free Gaussian distribution with bi-free characteristic triplet $(\mathbf{v},\mathbf{A},0)$ is non-full if and only if
$\mathbf{A}$ is singular, in which case, it is supported on the line $\langle\mathbf{u},(s,t)\rangle=0$, where $\mathbf{u}$ is a nonzero vector in the kernel of $\mathbf{A}$.
\end{cor}

\begin{cor} A bi-free compound Poisson distribution with rate $\lambda>0$ and jump distribution $\nu$ is non-full if and only if $\nu$ is
$\mathcal{M}_{\mathbb{R}^2}^\mathbf{0}$-nonfull, in which case they are supported on the same line. Consequently, any bi-free Poisson distribution is non-full.
\end{cor}

\begin{pf} Following the notations in (\ref{compound1}) and (\ref{compound2}), Theorem \ref{fullMfull2} yields that $\nu$ is supported on the line $\langle\mathbf{u},(s,t)\rangle=0$ for some $\mathbf{u}\in\mathbb{R}^2$ if and only if $\mu$ is supported on the line $\langle\mathbf{u},(s,t)\rangle=\langle\mathbf{u},\mathbf{v}\rangle=0$, as desired.
\end{pf} \qed

\begin{thm} Let $\mu\in\mathscr{P}_{\mathbb{R}^2}$ be $*$-infinitely divisible. Then $\mu$ is full if and only if $\Lambda(\mu)$ is full.
\end{thm}

\begin{pf} Recall that $\mu$ is non-full if and only if its characteristic function has the property that $|\widehat{\mu}(\lambda\mathbf{u})|=1$ for all $\lambda\in\mathbb{R}$, where $\mathbf{u}$ is some nonzero vector. If $P$ is the Poisson part of $\mu$, then  \[|\widehat{\mu}(\lambda\mathbf{u})|=|\widehat{P}(\lambda\mathbf{u})|\exp
\left[-\frac{1}{2}\lambda^2\langle\mathbf{Au,u}\rangle\right]\] yields that $\mu$ is non-full if and only if $\mathbf{Au}=0$ and $|\widehat{P}(\lambda\mathbf{u})|=1$ for all $\lambda$. Then the desired result follows from Theorem \ref{fullMfull2} and (\ref{bijection}).
\end{pf} \qed

\subsection*{Acknowledgement} The first-named author is supported by JSPS Grant-in-Aid for Young Scientists (B) 15K17549.
The second-named author was supported through a grant from the Ministry of Science and Technology in Taiwan MOST-104-2115-M-110-011-MY2 and a faculty startup grant from the National Sun Yat-sen University. The third-named author was supported by the NSERC Canada Discovery Grant RGPIN-2016-03796.

\end{document}